\documentclass[11pt, reqno]{amsart}
\usepackage[utf8]{inputenc}
\usepackage[T1]{fontenc}
\usepackage[english]{babel}

\usepackage[english]{d-omega3}

\usepackage{tikz}
\usetikzlibrary{decorations.pathreplacing}
\usetikzlibrary{arrows.meta}

\usepackage{mathtools}
\usepackage{amssymb}
\usepackage{graphicx}
\usepackage[all]{xy}

\usetikzlibrary{shapes.geometric}

\begin{document}

%%%%%%%%%%%%%%%% Insert metadata of the paper

\selectlanguage{english}
%\selectlanguage{ukrainian}

% Set the title of the paper:
% If yout title is long and will not fit to hte top margin of the paper
% then specify a short title in square brackets
\title[Cantorval property for a natural base]{A criterion for the set of subsums to be a Cantorval for a certain class of multigeometric series}

% Your title may contain line breaks like '\\'
% In that case define the following command including there your title without such symbols '\\'
% \contentsTitle{}

%%% information on the first author:
%%  \author[Name which will be printed at the end of the paper]
%%         {Name printed in the paper title}

\author[M. Pratsiovytyi]{Mykola Pratsiovytyi}
\organization{Institute of Mathematics of NAS of Ukraine, Dragomanov Ukrainian State University}
\email{prats4444@gmail.com}
\orcid{0000-0001-6130-9413}
%\thanks{This work was supported by a grant from the Simons Foundation (SFI-PD-Ukraine-00014586, M.P., O.V.).}
%%%%%%%%%%%% information on the second author:

\author[Ya. Vynnyshyn]{Yaroslav Vynnyshyn}
\organization{Kyiv Natural Science Lyceum No.~145}
\email{vynnyar@gmail.com}

\author[O. Vynnyshyn]{Oleh Vynnyshyn}
\organization{Institute of Mathematics of NAS of Ukraine}
\email{oleh.vynnyshyn@imath.kiev.ua}
\orcid{0009-0006-7113-1501}

%\shortAuthorsList{M.~Pratsiovytyi, O.~Vynnyshyn, Ya.O.~Vynnyshyn}
%%%%%%%%%%%% abstracts
\abstract{english}{The paper studies the topological structure of the set of subsums (incomplete sums) of multigeometric series of the form $\frac{k_1}{s} + \frac{k_2}{s} + \ldots + \frac{k_m}{s} + \frac{k_1}{s^2} + \frac{k_2}{s^2} + \ldots+ \frac{k_m}{s^2} + \frac{k_1}{s^3} + \frac{k_2}{s^3} +\ldots$, where $k_1, \ldots, k_m, s \in \mathbb{N}$.
In terms of graph theory, necessary and sufficient conditions are established for the set of subsums of a given series not to be a finite union of closed intervals, as well as necessary and sufficient conditions for this set to contain an interval. As a consequence of these statements, a criterion for the set of subsums to be a Cantorval is obtained.
An example of applying the criterion to a series whose Cantorval property was previously unknown is given.

 A more general object of study of topological properties is the set of numbers whose representations in a numerical system with a natural base and a redundant alphabet are subject to restrictions on the use of digits. It is for this object that the above statements are proved.
}

% \dedication{...}
\thanks{This work was supported by a grant from the Simons Foundation (SFI-PD-Ukraine-00014586, P.M., V.O.).}
\msc{40A05, 28A80, 11A67, 05C90}
%
% \doi{}
% \received{}
% \accepted{}

\maketitle

%%% Text of the paper

\section{Introduction}
Let $r$ be a fixed natural number, $A_r = \{0,1, \ldots, r\}$ an alphabet; $L_r = A_r \times A_r \times A_r \times \cdots$ the set of sequences of elements of the alphabet, in particular $L_1$ is the set of sequences of zeros and ones; $x_1 + x_2 + \cdots + x_k+ \cdots$ an absolutely convergent numerical series with real nonzero terms.

Recall~\cite{Kakeya_1914_partial_sums} that the set of subsums (incomplete sums) of a given series is defined as

\[ \mathcal{E}\bigl(x_n\bigr) \coloneqq
\left\{ x \; : \;
x = \sum_{n=1}^{\infty} c_n x_n, \quad
(c_n) \in L_1 \right\}. \]

The topological properties of sets of incomplete sums of absolutely convergent series were first studied by Soichi Kakeya~\cite{Kakeya_1914_partial_sums}, who proved that they are perfect; found necessary and sufficient conditions under which the set of subsums is an interval or a finite union of intervals, and also found sufficient conditions for the set of subsums to be nowhere dense. Later, in the paper~\cite{3_Guthrie_Nymann} a classification result was obtained:

\begin{theorem}[Guthrie--Nymann classification theorem]
\label{Hatri_Nymann_th}
The set of incomplete sums (subsums) of an absolutely convergent series is

\begin{enumerate}
%\item
%either a finite set
\item
either a union of a finite number of intervals
\item
or homeomorphic to the Cantor set
\item
or a Cantorval, that is, a set homeomorphic to the set
\[C \,\cup\, \bigcup_{n=1}^{\infty} G_{2n-1},\]
where $C$ is the classical Cantor set and $G_n$ are the intervals
removed at the $n$-th step of its construction.
\end{enumerate}
\end{theorem}

The notion of a Cantorval was introduced in the paper~\cite{Mendes_Oliveira_MCantorvals_1994}; the most well-known example is the Guthrie--Nymann Cantorval~\cite{3_Guthrie_Nymann}, which is the set of incomplete sums of the series with terms $x_{2n-1} = \frac{2}{4^n}, x_{2n} = \frac{3}{4^n}$.

Multigeometric series of the form
%\[ \sum_{n=1}^{\infty} \sum_{j=1}^{m} k_j q^n, \]
\[k_1 q + k_2 q + \ldots + k_m q + k_1 q^2 + k_2 q^2 + \ldots+ k_m q^2 + k_1 q^3 + \ldots,\]
where $ k_1,\ldots,k_m > 0, \quad q \in (0,1)$
have been studied in the papers
~\cite{Ferens1984,
Nitecki2015,
Bartoszewicz_Filipczak_Szymonik_2014,
Karvatskyi_Murillo_Viruel_2024_generalized_multigeom,
Pratsiovytyi_Karvatskyi_2023_trigonometric,
PratsiovytyiKarvatskyiMakarchuk2026}
and others. For them, in the paper~\cite{Bartoszewicz_Filipczak_Szymonik_2014} sufficient conditions for the set of incomplete sums to be a Cantorval are given, and the notation $(k_1, k_2, \ldots, k_m; q)$ for a multigeometric series is introduced, which we shall use hereafter.
In the paper~\cite{Nowakowski_OnANewCondition} a new sufficient condition for the set of subsums to be a Cantorval is obtained.

The current state and a list of open problems concerning the topological and metric properties of sets of incomplete sums of convergent series are presented in the survey~\cite{GlabPrusWisniowski2025}, where one of the central problems indicated is a criterion for the set of incomplete sums of a series to be a Cantorval.

The object of study in the present paper is
the series

\begin{equation}\label{multigeom_series}
(k_1, k_2, \ldots, k_m; \frac{1}{s}) = \frac{k_1}{s} + \ldots + \frac{k_m}{s} + \frac{k_1}{s^2} + \ldots+ \frac{k_m}{s^2} + \ldots,
\end{equation}
where $s,k_1, \ldots, k_m$ are fixed natural numbers, $s \geq 2$, as well as the sets associated with it.
Historically, it was precisely multigeometric series with integer coefficients $k_1, \ldots, k_m$ that provided examples of sets of subsums that are neither finite unions of intervals nor Cantor-type sets~\cite{3_Guthrie_Nymann,Ferens1984, Jones2011}.
Sets of subsums of multigeometric series with integer coefficients also arise naturally in probability theory~\cite{PratsiovytyiKarvatskyiMakarchuk2026, PratsiovytyiRatushniak2025}: they are the spectra of the corresponding infinite Bernoulli convolutions.

The criterion obtained in this paper for the Cantorval property for the class of multigeometric series (\ref{multigeom_series}) just described characterizes the topological type of their sets of subsums, and thus can be viewed as a partial solution to the problem of the topological structure of the set of subsums of a series.

\section{Sets of numbers with restrictions on the use of digits in their representations in a system with a natural base and a redundant alphabet}

%\subsection{Number systems with a redundant alphabet}

Henceforth $s$ and $r$ are fixed natural numbers, $2 \leq s \leq r$. It is known~\cite{Renyi1957} that for any number $x \in [0, \frac{r}{s-1}]$ there exists at least one sequence $(\alpha_k) \in L_r$ such that $x = \sum_{k=1}^{\infty}\alpha_ks^{-k } = \Delta_{\alpha_1 \alpha_2 \ldots \alpha_n \ldots}^{r_s}$.
The latter symbolic notation $\Delta_{\alpha_1 \alpha_2 \ldots \alpha_n \ldots}^{r_s}$
is called the $r_s$-representation of the number $x$, and $\alpha_n$ is the $n$-th digit of this representation. Since the alphabet $A_r$ is redundant, in general numbers of the interval $[0, \frac{r}{s-1}]$ may have a unique representation, a finite, a countable, or even a continuum of distinct representations~\cite{PratsiovytyiVynnyshyn2026, Pratsiovytyi_Ratushniak_2026}.

The main object of consideration in this paper is the set
\[C = C[r_s,V] = \left\{ x : x = \Delta_{\alpha_1 \alpha_2 \ldots \alpha_n \ldots}^{r_s}, \alpha_n \in V \subset A_r \right\},\]
where $V$ is a fixed proper non-singleton subset of the alphabet $A_r$. In this paper we are interested in the topological properties of the set $C$ under the condition $\{0,r\} \subset V$.

For a multigeometric series $(k_1, k_2, \ldots, k_m; \frac{1}{s})$, where $k_i, s$ are natural numbers, the set of incomplete sums $\mathcal{E}(x_n)$ is a set $C[r_s, V]$ provided that $V = \{v : v = \sum_{i=1}^m c_i k_i, \quad (c_1,c_2, \ldots , c_m) \in A_1^m \}$.

\emph{Example 1}. The series $(3,2; \tfrac{1}{4})$ (Guthrie--Nymann~\cite{3_Guthrie_Nymann}). For it, $V = \{0,2,3,5\}$, and the set of incomplete sums $C$ is a Cantorval.

\emph{Example 2}. The series $(6,1; \tfrac{1}{4})$ (Kenyon's example~\cite{Nitecki2015}). In this case $V = \{0,1,6,7\}$, and $C$ is a Cantorval.

\emph{Example 3}. The series $(1,8; \tfrac{1}{4})$~\cite{GlabPrusWisniowski2025}. The alphabet $V = \{0,1,8,9\}$, $C = [0,1] \cup [2,3]$.

Not every set $C[r_s,V]$ is the set of subsums of some multigeometric series of the form \eqref{multigeom_series}. Also, in general, the Guthrie--Nymann classification theorem does not apply to such sets.

\section{Structure of the set $C[r_s,V]$}

Henceforth we shall use the notation $\nabla_{\alpha_1 \alpha_2 \ldots \alpha_m}^s $ for the $s$-adic cylindrical half-interval of rank $m$ with base $\alpha_1 \alpha_2 \ldots \alpha_m$, that is

\[ \nabla_{\alpha_1 \alpha_2 \ldots \alpha_m}^s :=
\left[ \sum_{k=1}^{m} \frac{\alpha_k}{s^k}, \;
\sum_{k=1}^{m} \frac{\alpha_k}{s^k} + \frac{1}{s^m} \right), \]
where $\alpha_k \in \{0,1,\ldots,s-1\}$.

Clearly,
\[ \nabla_{\alpha_1 \alpha_2 \ldots \alpha_m}^s =
\nabla_{\alpha_1 \alpha_2 \ldots \alpha_m0}^s \cup
\nabla_{\alpha_1 \alpha_2 \ldots \alpha_m1}^s \cup \ldots
\cup \nabla_{\alpha_1 \alpha_2 \ldots \alpha_m[s-1]}^s
\]
%for a cylindrical half-interval of rank $m$ with base $\alpha_1 \alpha_2 \ldots \alpha_m$.

The half-open interval $[0;1)$ will be called the cylindrical half-interval of rank zero. It is accordingly partitioned into $s^n$ pairwise disjoint cylindrical half-intervals of rank $n$ of length $\frac{1}{s^n}$.

For convenience we denote
\[ {\Delta}_{\alpha_1 \alpha_2 \ldots \alpha_n} \coloneqq
C[r_s,V] \cap
\nabla_{\alpha_1 \alpha_2 \ldots \alpha_n}^s.\]
%We shall denote by $A \oplus B \coloneqq \{a + b \mid a \in A,\ b \in B\}$ the arithmetic sum of the sets $A$ and $B$.
Also let $[0,M)$ be the smallest half-interval containing the set $C[r_s, V]$, taking into account that $0 \in V$. It is clear that $M = [\frac{r}{s-1}] + 1$.

If we introduce the notation for the sets
\[ E_n \coloneqq \left\{ x \; : \;
x = \sum_{k=1}^{n} \frac{\alpha_k}{s^k}, \quad
(\alpha_1, \ldots, \alpha_n) \in V^n \right\}, \]

\[ R_n \coloneqq \left\{ x \; : \;
x = \sum_{k=n+1}^{\infty} \frac{\alpha_k}{s^k}, \quad
\alpha_k \in V, \forall k \in \mathbb{N} \right\}, \]

then clearly $E_n \subset E_{n+1} \subset C, \, R_{n+1} \subset R_n \subset C \quad \forall n \in \mathbb{N}$, in particular
\[C = E_n \oplus R_n, \quad C = [\bigcup_{k=1}^{\infty} E_k] = [\lim_{k \to \infty} E_k],  \]
where $[B]$ denotes the closure of the set $B$.

The elements of the set $E_n$ will be called \emph{points of rank $n$}.
\begin{lemma}\label{points_rank_n}
 Each point of rank $n$ can be represented in the form $\frac{a}{s^n}$, where $a$ is some nonnegative integer.
\end{lemma}
\begin{proof}
Indeed, each element of the set $E_n$ has the form
\[\frac{\alpha_1}{s} + \frac{\alpha_2}{s^2} + \ldots + \frac{\alpha_n}{s^n} = \frac{\sum_{i=1}^n \alpha_is^{n-i}}{s^n},\]
where $\alpha_i \in V$ are integers.
\end{proof}

The set of points of rank $0$ will be taken to be the singleton $\{0\}$.

\medskip

Consider the $M$-point set $T$ associated with the cylindrical half-interval $\nabla_{\alpha_1 \alpha_2 \ldots \alpha_n}^s$, which contains the numbers defined by the equality

\[t_j := \frac{\alpha_1}{s} + \ldots + \frac{\alpha_n}{s^n} - \frac{M-1}{s^n} + \frac{j}{s^n}, \quad j \in \{0,1,\ldots, M-1\}.\]

\medskip
\begin{tikzpicture}[
  scale=1.05,
  axis/.style={->, line width=1pt},
  dot/.style={circle, fill, inner sep=1.5pt},
  hollow/.style={circle, draw, inner sep=2pt},
  cyl/.style={line width=2pt, black}
]

% axis
\draw[axis] (0,0) -- (8,0);

% cylinder
\draw[cyl] (6,0.12) -- (7,0.12);
\draw[dashed] (6,0) -- (6,0.12);
\draw[dashed] (7,0) -- (7,0.12);
\node[above] at (6.5,0.18) {$\nabla_{\alpha_1 \alpha_2 \ldots \alpha_n}^s$};

% prefix points
\node[dot] at (0,0) {};
\node[dot] at (1,0) {};
%\node[hollow] at (2,0) {};
\node[dot] at (2,0) {};
\node[dot] at (4,0) {};
%\node[hollow] at (5,0) {};
\node[dot] at (5,0) {};
\node[dot] at (6,0) {};

% labels
\node[below] at (0,-0.1) {\scriptsize $t_0$};
\node[below] at (1,-0.1) {\scriptsize $t_1$};
\node[below] at (2,-0.1) {\scriptsize $t_2$};

\node[below] at (4,-0.1) {\scriptsize $t_{M-3}$};
\node[below] at (5,-0.1) {\scriptsize $t_{M-2}$};
\node[below] at (6,-0.1) {\scriptsize $t_{M-1}$};

\node at (3,0) {$\cdots$};

% distance 1/s^n
\draw[
  decorate,
  decoration={brace, amplitude=6pt}
] (1,0.1) -- (2,0.1)
node[above] at (1.5,0.2) {$\frac{1}{s^n}$};

\end{tikzpicture}

The point $t_{M-1}$ is the left endpoint of the cylindrical half-interval $\nabla_{\alpha_1 \alpha_2 \ldots \alpha_n}^s$, and the points $t_{j-1}$ are obtained by shifting the points $t_j$ by $\frac{1}{s^n}$ to the left.

\begin{lemma}\label{use_only_points_t_i}
For an arbitrary $\nabla_{\alpha_1 \alpha_2 \ldots \alpha_n}^s$ the equality

\[
\Delta_{\alpha_1 \alpha_2 \ldots \alpha_n}
=
\bigl( (E_n \cap T)
\oplus R_n \bigr)
\cap \nabla_{\alpha_1 \alpha_2 \ldots \alpha_n}^s.
\]
holds.
\end{lemma}
\begin{proof}
By definition
\[\Delta_{\alpha_1 \alpha_2 \ldots \alpha_n} =
C \cap \nabla_{\alpha_1 \alpha_2 \ldots \alpha_n}^s =
\bigl( E_n \oplus R_n \bigr) \cap
[t_{M-1}, t_{M-1} + \tfrac{1}{s^n}). \]
Note that for arbitrary $A \subset \mathbb{R}$, $B \subseteq [b_1, b_2)$
\[\bigl( A \oplus B\bigr) \cap [x,y) =
\Bigl( \bigl(A \cap (x-b_2,\, y - b_1)\bigr)\oplus B \Bigr)
\cap [x,y). \]
Since the set $R_n \subset[0,\frac{M}{s^n})$, we have
\[\Delta_{\alpha_1 \alpha_2 \ldots \alpha_n} =
\Bigl(
    \bigl(
        E_n \cap (t_{M-1} - \tfrac{M}{s^n},\,
        t_{M-1} + \tfrac{1}{s^n})
    \bigr)
    \oplus R_n
\Bigr)
\cap [t_{M-1}, t_{M-1} + \tfrac{1}{s^n}).\]
 By the definition $t_0 - \tfrac{1}{s^n} = t_{M-1} - \frac{M}{s^n}$, we have
\[
\Delta_{\alpha_1 \alpha_2 \ldots \alpha_n} = \bigl( (E_n \cap(t_0 - \tfrac{1}{s^n},
t_{M-1} + \tfrac{1}{s^n}))
\oplus R_n \bigr)
\cap \nabla_{\alpha_1 \alpha_2 \ldots \alpha_n}^s.\]

Each of the points $t_0, \ldots, t_{M-1}$ may or may not belong to the set $E_n$, while no other point from $(t_0 - \tfrac{1}{s^n}, t_{M-1} + \tfrac{1}{s^n})$ belongs to $E_n$ (Lemma~\ref{points_rank_n}). Therefore
\[\Delta_{\alpha_1 \alpha_2 \ldots \alpha_n} =
\bigl( (E_n \cap \{t_0, t_1, \ldots, t_{M-1}\})
\oplus R_n \bigr)
\cap \nabla_{\alpha_1 \alpha_2 \ldots \alpha_n}^s.
\qedhere \]
\end{proof}

\section{Sets of prefix points and prefixes}

\begin{definition}
The \emph{set of prefix points} of a cylindrical half-interval
$\nabla_{\alpha_1 \ldots \alpha_n}^s$ is the set
\[ B(\nabla_{\alpha_1 \ldots \alpha_n}^s) :=
\{t_0, t_1, \ldots, t_{M-1}\} \cap E_n. \]
\end{definition}

\begin{definition}
The \emph{prefix} of a cylindrical half-interval $\nabla_{\alpha_1\ldots\alpha_n}^s$ of rank $n$ is the ordered tuple $(p_0, p_1, \ldots, p_{M-1})$, where
$p_i=1$ if $t_i\in E_n$, and $p_i=0$ if $t_i\notin E_n$.
\end{definition}

%That is, the set of prefix points determines the prefix:
%\begin{tikzpicture}[
%  scale=1.05,
%  axis/.style={->, line width=1pt},
%  dot/.style={circle, fill, inner sep=1.5pt},
%  hollow/.style={circle, draw, inner sep=2pt},
%  cyl/.style={line width=2pt, black}
%]

%% axis
%\draw[axis] (0,0) -- (8,0);

%% cylinder
%\draw[cyl] (6,0.12) -- (7,0.12);
%\draw[dashed] (6,0) -- (6,0.12);
%\draw[dashed] (7,0) -- (7,0.12);

% prefix points
%\node[dot] at (0,0) {};
%\node[dot] at (1,0) {};
%\node[hollow] at (2,0) {};

%\node[dot] at (4,0) {};
%\node[hollow] at (5,0) {};
%\node[dot] at (6,0) {};

%texts
%\node[below] at (0,-0.1) {\scriptsize $t_0$};
%\node[below] at (1,-0.1) {\scriptsize $t_1$};
%\node[below] at (2,-0.1) {\scriptsize $t_2$};

%\node[below] at (4,-0.1) {\scriptsize $t_{M-3}$};
%\node[below] at (5,-0.1) {\scriptsize $t_{M-2}$};
%\node[below] at (6,-0.1) {\scriptsize $t_{M-1}$};

%\node[above] at (0,0.1) {\scriptsize $1$};
%\node[above] at (1,0.1) {\scriptsize $1$};
%\node[above] at (2,0.1) {\scriptsize $0$};

%\node[above] at (4,0.1) {\scriptsize $1$};
%\node[above] at (5,0.1) {\scriptsize $0$};
%\node[above] at (6,0.1) {\scriptsize $1$};
%\end{tikzpicture}

The set of all possible prefixes is denoted $ P := \{0,1\}^M$. We note that there are finitely many prefixes. Therefore, for a fixed rank $n$ there exist only finitely many pairwise distinct (up to a parallel translation) sets $\Delta_{\alpha_1\ldots\alpha_n}$, since for the same prefix the sets are equal up to a parallel translation.
For different ranks $n$ and $m$ --- if $\nabla_{\alpha_1\ldots\alpha_n}^s$ and
$\nabla_{\beta_1\ldots\beta_m}^s$
have the same prefixes, then
$\Delta_{\alpha_1\ldots\alpha_n}$ and $\Delta_{\beta_1\ldots\beta_m}$
are similar with coefficient $s^{m-n}$.

\subsection*{Transition from rank $n$ to rank $n+1$}

We describe the transition from cylindrical half-intervals of rank $n$ (in the representation $C = E_n \oplus R_n$) to cylindrical half-intervals of rank $n+1$ (in the representation $C = E_{n+1} \oplus R_{n+1}$).

%$l = \sum_{i=1}^{n} \frac{\alpha_i}{s^i} =\Delta_{\alpha_1 \ldots \alpha_n (0)}$
Let $(\alpha_1 \alpha_2 \ldots \alpha_n)$ be a fixed tuple of elements of the alphabet $A_s$, $l = \sum_{i=1}^{n} \frac{\alpha_i}{s^i}$. Then for an arbitrary $d \in A_s$ we have

\begin{equation}\label{transition_B}
B(\nabla_{\alpha_1 \alpha_2 \ldots \alpha_n d}^s) =
\Bigl(B(\nabla_{\alpha_1 \alpha_2 \ldots \alpha_n}^s)
\oplus \tfrac{V}{s^{n+1}} \Bigr)
\cap \Bigl[l +\tfrac{d-M+1}{s^{n+1}},l + \tfrac{d}{s^{n+1}}\Bigr]
\end{equation}

\begin{equation}\label{transition_SetOfSubsums}
\Delta_{\alpha_1 \ldots \alpha_n d} =
\Bigl( B(\nabla_{\alpha_1 \ldots \alpha_n d}^s) \oplus R_{n+1} \Bigr)
\cap \nabla_{\alpha_1 \ldots \alpha_n d}^s.
\end{equation}

The prefix points $\{t_0,\ldots,t_{M-1}\}$ of the subcylinder $\nabla_{\alpha_1 \ldots \alpha_n d}$ lie in the interval $[l +\tfrac{d-M+1}{s^{n+1}},l + \tfrac{d}{s^{n+1}}]$, both endpoints being themselves prefix points ($t_0$ and $t_{M-1}$).
The prefixes of subcylinders (direct descendants) are uniquely determined by the prefix of their direct ancestor. This result is natural, since the prefix of a cylinder completely determines the structure of $C$ on it, and hence on all of its subcylinders.

We note that the above reasoning about the structure of the set $C$ remains valid if the initial (rank $0$) cylindrical half-interval is taken to be not only $[0,1)$, but any half-interval $[k,k+1)$, $k \in \{0,1,\ldots,M-1\}$; the scheme described can be used to study the topological structure of the set $C$ on each of them.

\section{Transition graphs}

\begin{definition}
Let $W$ be a nonempty set. A relation $R \subset W \times W$ is called a \emph{partial equivalence relation}~\cite{Mitchell_1996_PER} if it is symmetric $(aRb \;\Rightarrow\; bRa)$ and transitive $(aRb \ \text{and}\ bRc \;\Rightarrow\; aRc)$.

%that is, for all $a,b,c \in W$ the following conditions hold:

\end{definition}
%Here reflexivity is not required in general, that is, $aRa$ may fail for some $a$.

For an element $a$ with $aRa$, the \emph{partial equivalence class} is the set ${b \in W : aRb}$.
In particular, we note that the set $W$ may contain elements without the property $aRa$ (reflexivity), that is, elements that do not belong to any partial equivalence class induced by $R$.

%\subsection*{Partial equivalence relation for a directed graph.}

Let $G$ be a directed graph. For vertices $a,b$ of the graph $G$ we set
$ a \sim b$,
if there exist directed paths from vertex $a$ to vertex $b$ and from vertex $b$ to vertex $a$. In other words, $a \sim b$ if there exists a cycle (a closed directed path) that passes through both $a$ and $b$.

The relation $\sim$ is a partial equivalence relation on the set of vertices of the graph, and its classes correspond to the nontrivial strongly connected components of the graph. Vertices that do not lie on any cycle do not belong to any partial equivalence class. We call them \emph{transient} vertices.

\begin{definition}\label{graph_of_transitions}
The \emph{transition graph} $G$ is the directed graph with vertex set $P = \{0,1\}^M$ (the vertices are all possible prefixes), in which from each vertex $q \in P$, for each $d \in \{0,1,\ldots,s-1\}$, there is an edge $q \rightarrowtail q_d$, where $q_d$ is the prefix of the $d$-th subcylinder of an arbitrary cylindrical half-interval with prefix $q$.
\end{definition}

Moving along an edge of the directed graph corresponds to choosing the corresponding subcylinder, so movement along the edges of the graph determines a sequence of nested subcylinders. The symbol $\rightarrowtail$ here and below denotes a single directed edge. The definition is well-posed: the prefix $q_d$ does not depend on the choice of the half-interval with prefix $q$ (see formula~\eqref{transition_B}).

From each vertex of the graph $G$ there emanate exactly $s$ edges; loops are allowed (the prefix of the $d$-th subcylinder may coincide with the prefix of the half-interval itself) as well as multiple edges (different values of $d$ may yield the same descendant prefix).

We shall write $a \rightsquigarrow b$ if in the graph $G$ there exists a directed path from vertex $a$ to vertex $b$.
As described above, the relation $\sim$ ($a \sim b$ if $a \rightsquigarrow b$ and $b \rightsquigarrow a$) partitions the vertices of the graph into partial equivalence classes and transient vertices. We contract each partial equivalence class (that is, each nontrivial strongly connected component of the graph) to a single vertex; this procedure is called the condensation of the graph \cite{cormen2022introduction}. As a result we obtain a directed acyclic graph $\bar{G}$, whose vertices are the partial equivalence classes and the transient vertices.

\begin{definition}
A partial equivalence class is called \emph{terminal} if all edges emanating from its vertices lead to vertices of the same class (in other words, no edge emanates from the corresponding vertex of the condensed graph $\bar{G}$).
\end{definition}

\begin{example}
The class consisting of the vertex \((0,0,\ldots,0)\) is terminal.
Indeed, the prefix $(0,\ldots,0)$ means that $B(\nabla_{\alpha_1\ldots\alpha_n}^s)= \varnothing$.
By~\eqref{transition_B} $B(\nabla_{\alpha_1\ldots\alpha_n d}^s) = \varnothing$ for all $d \in \{0,\ldots,s-1\}$, that is, the prefixes of the subcylinders are also $(0,\ldots,0)$. Hence all $s$ edges from this vertex are loops: the vertex $(0,\ldots,0)$ forms a class, and no edge leaves its bounds.
\end{example}

\begin{definition}
The \emph{level} of a partial equivalence class or a transient vertex is the maximum length of a directed path in the condensed graph $\bar{G}$ that starts at the given class/vertex and ends at some terminal class.
\end{definition}

\textbf{Remark.} Terminal classes are classes of level $0$.

\textbf{Remark.} If from a class/transient vertex $X$ one can pass to another class/transient vertex $Y$, then the level of $X$ is strictly greater than the level of $Y$.

\begin{definition}
\emph{Starting prefixes} are the prefixes of the cylindrical half-intervals of rank zero $[k-1, k)$, $k = 1, \ldots, M$, that is, the tuples
\[
(0,0,\ldots,0,1),\ (0,0,\ldots,1,0),
\ \ldots,\ (1,0,\ldots,0,0).
\]
In the graph $G$ the starting prefixes will be called \emph{starting vertices}.
\end{definition}

By the construction of the graph $G$, for each $n$ the prefixes of the cylindrical half-intervals of rank $n$ contained in $[k-1, k)$ are precisely those vertices to which there exists a path of length $n$ in $G$ from the starting vertex of the half-interval $[k-1, k)$. Therefore, to study the topological structure of the set $C$ it suffices to consider the starting vertices and the vertices reachable from them.

\begin{definition}
By $G'$ we denote the subgraph of $G$ with vertex set $P' \subseteq P$ consisting of the starting vertices and the vertices reachable from them by directed paths. The edges of $G'$ are all edges of $G$ that emanate from vertices of $P'$; by $\bar{G'}$ we denote the condensation of the graph $G'$.

Together with each vertex, the graph $G'$ contains all of its descendants, therefore the terminality of classes and the levels of vertices in $\bar{G'}$ are the same as in $\bar{G}$. Henceforth we work with the graphs $G'$ and $\bar{G'}$.
\end{definition}

\section{A criterion for the presence of an interval in $C$}

\begin{definition}
A \emph{gap} of the set $C$ is an open interval
$
(a,b)\subset [0,M], % possibly replace with [0, r/(s-1)]
\qquad
(a,b)\cap C=\varnothing.
$
\end{definition}
\begin{remark}
In some papers a gap is defined as a maximal open interval in the complement of $C$~\cite{Nowakowski_OnANewCondition}. We use the term in a broader sense: ``any open interval that does not intersect $C$''.
\end{remark}

\begin{lemma}\label{lemma_gap}
Let $\nabla_{\alpha_1\ldots\alpha_n}^s \subset [0,\frac{r}{s-1}]$ have prefix $q$. The half-interval $\nabla_{\alpha_1\ldots\alpha_n}^s$ contains a gap if and only if in the prefix graph $G'$ there exists a directed path $ q \;\rightsquigarrow\; (0,0,\ldots,0).$
\end{lemma}

\begin{proof}

\emph{Sufficiency $(\Leftarrow)$.}
Let the path $q \rightsquigarrow (0,\ldots,0)$ have length $k$. By the construction of the graph this path corresponds to a choice of digits $d_1,\ldots,d_k$, so the cylindrical subinterval $\nabla_{\alpha_1\ldots\alpha_n d_1\ldots d_k}^s$ has prefix $(0,\ldots,0)$. Then $B(\nabla_{\alpha_1\ldots\alpha_n d_1\ldots d_k}^s) = \varnothing$ and, accordingly, $\Delta_{\alpha_1\ldots\alpha_n d_1\ldots d_k} = \varnothing$.

\emph{Necessity $(\Rightarrow)$.} Let $\nabla_{\alpha_1\ldots\alpha_n}^s$ contain a gap $(a,b)$. Choose $N > n$ so that the inequality $\tfrac{1}{s^N} < \tfrac{b-a}{3}$ holds. For such $N$ there exists a subcylinder
$\nabla_{\alpha_1\ldots\alpha_n\beta_{n+1}\ldots\beta_N}^s \subset (a,b)$.
Choose a number $k$ large enough that $s^k> M+1$. Then iteratively $k$ times we choose the rightmost subcylinder in $\nabla_{\alpha_1\ldots\alpha_n\beta_{n+1}\ldots\beta_N}^s$, that is, $\nabla_{\alpha_1\ldots\alpha_n\beta_{n+1}\ldots\beta_N [s-1]\ldots[s-1]}^s$. It and all $M$ of its prefix points lie in $\nabla_{\alpha_1\ldots\alpha_n\beta_{n+1}\ldots\beta_N}^s$, and hence lie in $(a,b)$. Therefore its prefix is $(0,0,\ldots,0)$. Thus we have indicated a directed path from $q$ to $(0,0,\ldots,0)$ in the graph $G'$.
\end{proof}

\begin{lemma}\label{lemma_black}
Let the cylindrical interval
$\nabla_{\alpha_1\ldots\alpha_n}^s$
have prefix $q$.
The half-interval $ \nabla_{\alpha_1\ldots\alpha_n}^s \subset C$
if and only if
in the prefix graph $G'$ there is no directed path
$
q \rightsquigarrow (0,0,\ldots,0).
$
\end{lemma}

\begin{proof}

\emph{Necessity $(\Rightarrow)$.} The existence of a path $q \rightsquigarrow (0,\ldots,0)$ implies the existence of a gap in $\nabla_{\alpha_1\ldots\alpha_n}^s$ (Lemma~\ref{lemma_gap}).

\emph{Sufficiency $(\Leftarrow)$.}
Suppose that there is no directed path $q \rightsquigarrow (0,0,\ldots,0)$.
Take an arbitrary point $x_0 \in \nabla_{\alpha_1\ldots\alpha_n}^s$ and consider the sequence of cylinders
$ \nabla_{\alpha_1\ldots\alpha_n }^s \supset
 \nabla_{\alpha_1\ldots\alpha_n \alpha_{n+1}}^s   \supset
 \nabla_{\alpha_1\ldots\alpha_n
\alpha_{n+1} \alpha_{n+2}}^s \supset \ldots$,
where $\alpha_{n+1},\alpha_{n+2},\ldots$ are the subsequent digits of some $s$-adic expansion of $x_0$.

%If $x_0 = \tfrac{r}{s-1}$, then $x_0 \in C$ (attained when all $\alpha_i = r$), so we assume henceforth that $x_0 < \tfrac{r}{s-1}$.
%If $x_0$ is $s$-adic rational and $x_0 > 0$, that is, has the form $\frac{a}{s^j}, a > 0$, then it has two distinct representations, namely with periods $(0)$ and $(s-1)$. For uniqueness we choose the representation with period $(0)$.

The nested cylinders form a sequence of nested half-intervals $[a_k,b_k)$. If some $[a_k,b_k)$ had prefix $(0,0,\ldots,0)$, then there would exist a path $q \rightsquigarrow (0,\ldots,0)$, which contradicts the assumption. Accordingly, they have nonzero prefixes. Hence each of the intervals $[a_k - \tfrac{M-1}{s^{n+k}}, a_k]$ contains at least one prefix point
$x_k \in C$. The inequality
$|x_0 - x_k|
\leq \tfrac{M-1}{s^{n+k}} + \tfrac{1}{s^{n+k}} =
\tfrac{M}{s^{n+k}}$
holds,
that is, $x_k \to x_0$. The set $C$ is closed (the scheme of the proof is in~\cite{Kakeya_1914_partial_sums}), therefore also $x_0 \in C$. The point $x_0$ was chosen arbitrarily from $\nabla_{\alpha_1\ldots\alpha_n}^s$, hence $\nabla_{\alpha_1\ldots\alpha_n}^s \subset C$.

\end{proof}

\begin{theorem}\label{theorem_segment_existence}
The set $C$ contains an interval if and only if the condensed graph $\bar{G'}$ contains a terminal partial equivalence class different from the class consisting of the vertex
$(0,0,\ldots,0)$.
\end{theorem}
\begin{proof}
\emph{Necessity $(\Rightarrow)$.} Let $C$ contain an interval $[c,d]$. Choose $n$ large enough that
 $\tfrac{1}{s^n} < \tfrac{d-c}{3}$.
 Then there exists $\nabla_{\alpha_1\ldots\alpha_n}^s \subset [c, d]$; denote its prefix by $q$. Since $\nabla_{\alpha_1\ldots\alpha_n}^s \subset C$, it contains no gaps, so by Lemma~\ref{lemma_gap} there is no path $q \rightsquigarrow (0,\ldots,0)$ in the graph $G'$.
%Since $q$ is the prefix of a cylinder, we have $q \in P'$.
Consider the class/transient vertex in the graph $\bar{G'}$ that contains $q$. Since in the condensed graph $\bar{G'}$ the level always decreases under transitions, from $q$ there exists a path to a terminal class which, as shown, is not $\{(0,\ldots,0)\}$.

\emph{Sufficiency $(\Leftarrow)$.} Let $T$ be a terminal class in $\bar{G'}$ different from $\{(0,\ldots,0)\}$, and let $q \in T$ be an arbitrary vertex. Since $T$ is terminal, there is no path $q \rightsquigarrow (0,0,\ldots,0)$. Since $q \in P'$, there exists some cylindrical half-interval with prefix $q$ and, according to the previous Lemma~\ref{lemma_black}, it lies entirely in $C$. Inside the half-interval one can choose an interval.
\end{proof}

%--------------------------------------------
\section{A criterion for the finiteness of the number of intervals in $C$}

\begin{definition}
We define the \emph{colors of partial equivalence classes and transient vertices} as follows:

\emph{Level $0$.}
The terminal class $\{(0,0,\ldots,0)\}$ is assigned the \emph{white} color, and all other terminal classes the \emph{black} color.

\emph{Levels $\geq 1$.}
A class (transient vertex) is assigned a color depending on the colors of its immediate descendants, namely:
\begin{itemize}
  \item \emph{white}, if all immediate descendants have white color;
  \item \emph{black}, if all immediate descendants have black color;
  \item \emph{black-and-white} in the remaining cases, that is, if among the immediate descendants there are both a black and a white one, or at least one black-and-white descendant.
\end{itemize}
\end{definition}
\begin{definition}
    The \emph{color of a vertex} from some equivalence class is defined to be the color of the equivalence class to which it belongs.
\end{definition}

\begin{remark}
    Intuitively: white color indicates the presence of gaps in the descendants, while black indicates the presence of intervals.
\end{remark}

\begin{remark}\label{remark_black_is_full}
A cylindrical interval $\nabla_X^s$ lies entirely in $C$ if and only if the prefix $X$ is black (Lemma~\ref{lemma_black}).
\end{remark}

%For a prefix (vertex) $X$ we shall denote by $\nabla_X^s$ an arbitrary cylindrical half-interval with prefix $X$. Since the set $\Delta_{\alpha_1\ldots\alpha_n}$ is determined, up to similarity, by the prefix alone, all properties of interest to us do not depend on the choice of such a cylinder.

\begin{remark}\label{remark_white_no_segment}
A cylindrical interval $\nabla_W^s$ with a white prefix $W$ contains no interval
from $C$.

%Indeed, if there existed an interval $[a,b] \subset \nabla_W^s \cap C$, then inside $[a,b]$ there would lie some cylindrical interval $\nabla_{q^*}^s$ of sufficiently large rank. Then $\nabla_{q^*}^s \subset C$. By the previous Remark~\ref{remark_black_is_full} $q^*$ is black, so there exists a path from a white $w$ to a black $q^*$, a contradiction.

We also note that the white color of a prefix indicates only the absence of intervals in the cylindrical half-interval, while at the same time it may contain many points of $C$.
\end{remark}

\begin{lemma}\label{lemma_color_inheritance}
The color of any vertex contains the colors of all of its descendants in the following sense:
if in the graph $G'$ there is a path between vertices $q_1 \rightsquigarrow q_2$, then
\begin{itemize}
  \item if $q_2$ is black, then $q_1$ is black or black-and-white;
  \item if $q_2$ is white, then $q_1$ is white or black-and-white;
  \item if $q_2$ is black-and-white, then $q_1$ is black-and-white.
\end{itemize}
\end{lemma}

\begin{proof}
Let $q'$ be the penultimate vertex of the path $q_1 \rightsquigarrow q_2$.
In $G'$ there is an edge $q' \rightarrowtail q_2$, so either $q'$ and $q_2$ lie in the same equivalence class and therefore have the same color, or in $\bar{G'}$ the class/transient vertex containing $q_2$
is an immediate descendant of the class/transient vertex containing $q'$, and then, by the definition of color, the color of vertex $q'$ contains the color of vertex $q_2$. By induction on the length of the path we obtain that the color of vertex $q_1$ contains the color of vertex $q_2$.
\end{proof}

\subsection{Global structure of black-and-white classes}

\begin{definition}
A partial equivalence class is called
\emph{cyclic} if from each of its vertices there emanates exactly one edge leading to a vertex of the same class (the remaining $s-1$ edges leave the class).
\end{definition}

\begin{lemma}\label{infinite_paths}
If $\gamma_i:a \rightsquigarrow b, \, a,b \in G', \, i \in \mathbb{N}$ are such paths that any pair differs at some step and the vertex $b$ is black-and-white, then the cylindrical half-interval $\nabla_{a}$ (with prefix $a$) contains infinitely many intervals from $C$.
\end{lemma}

\begin{proof}
Since $b$ is black-and-white, from it there exists a path to the terminal class $\{(0,\ldots,0)\}$ (a white descendant) and a path to another terminal class (a black descendant). By Lemmas~\ref{lemma_gap} and~\ref{lemma_black}, a cylindrical interval with prefix $b$ contains both an interval (in $C$) and a gap.

The paths $\gamma_i$ differ at some step --- none is a prefix of another --- therefore the descendant cylinders of the half-interval $\nabla_{a}$ obtained by the transitions $\gamma_i$ are disjoint. For $\gamma_1, \ldots, \gamma_k$ we arrange the descendant cylinders on the real line, and in the cylinders in odd positions we select an interval from $C$, and in even positions a gap. Since $k$ is arbitrary, $\nabla_{a}$ contains infinitely many intervals from $C$.
\end{proof}

\begin{lemma}\label{necessary_condition1_branching}
If in the graph $G'$ there exists a non-cyclic black-and-white class,
then the set $C$ contains infinitely many intervals.
\end{lemma}

\begin{proof}
Let $K \in \bar{G'}$ be a non-cyclic black-and-white class, $x \in K$ an arbitrary vertex.
%and a cylindrical half-interval $\nabla_x$.
By the definition of non-cyclicity there exists a vertex $y \in K$ that has two internal outgoing edges $y \rightarrowtail a_1$, $y \rightarrowtail a_2$, where $a_1, a_2 \in K$. The class $K$ is a strongly connected component, therefore there exist paths $x \rightsquigarrow y$, $a_1 \rightsquigarrow x$ and $a_2 \rightsquigarrow x$, and we choose them so that they have the vertex $x$ only at the endpoints. We compose them as follows:
\[\gamma_{a_1} = x \rightsquigarrow y \rightarrowtail a_1 \rightsquigarrow x; \quad
\gamma_{a_2} = x \rightsquigarrow y \rightarrowtail a_2 \rightsquigarrow x.
\]

\begin{center}
\begin{tikzpicture}[
  dot/.style={circle, fill, inner sep=2pt},
  cls/.style={draw, line width=1.2pt},
  arr/.style={->, line width=1.1pt}
]

\draw[cls] (0,0) ellipse (2.6 and 1.4);
\node at (3.5,0) {$K$};

\node[dot] (x1) at (-1.5,0) {};
\node[left] at (x1) {$x$};
\node[dot] (x2) at (0,-0) {};
\node[dot] (x3) at (1.5,0) {};
\node[right] at (x3) {$y$};

\node[dot] (x41) at (1,1) {};
\node[right] at (x41) {$a_1$};
\node[dot] (x42) at (0.75,-1) {};
\node[right] at (x42) {$a_2$};
\node[dot] (x52) at (-0.75,-1) {};

\draw[arr, blue] (x1) to[bend left=15] (x2);
\draw[arr, red]  (x1) to[bend right=15] (x2);
\draw[arr, blue] (x2) to[bend left=15] (x3);
\draw[arr, red]  (x2) to[bend right=15] (x3);

\draw[arr,blue] (x3) -- (x41);
\draw[arr,blue] (x41) -- (x1);

\draw[arr,red] (x3) -- (x42);
\draw[arr,red] (x42) -- (x52);
\draw[arr,red] (x52) -- (x1);

\node[blue] at (-0.1,0.8) {$\gamma_{1}$};
\node[red] at (-0.1,-0.8) {$\gamma_{2}$};

\end{tikzpicture}
\end{center}
We set $\gamma_k = (\gamma_{a_1})^k \gamma_{a_2} : x \rightsquigarrow x$. Any pair differs at some step, the vertex $x$ is black-and-white, therefore by Lemma~\ref{infinite_paths} $\nabla_{\alpha_1\ldots\alpha_n}^s$ contains infinitely many intervals from $C$.
\end{proof}

\begin{lemma}\label{necessary_condition2_connecting}
If for two distinct black-and-white classes from $G'$ there exists a directed path from one of them to the other, then the set $C$ contains infinitely many intervals.
\end{lemma}
\begin{proof}
Let for distinct black-and-white classes $X, Y$ there exist a path $\gamma_{xy} \colon x \rightsquigarrow y$, where $x \in X$, $y \in Y$. Denote by $\gamma'_{xy}$ some subpath of it that contains $x$ and $y$ only at the endpoints. Inside the class $X$ we fix some path $\gamma_{xx} \colon x \rightsquigarrow x$ that has $x$ only at the endpoints.

\begin{center}
\begin{tikzpicture}[
  dot/.style={circle, fill, inner sep=2pt},
  cls/.style={draw, line width=1.2pt},
  arr/.style={->, line width=1.1pt}
]

\draw[cls] (0,0) ellipse (1.6 and 0.9);
\node at (2.5,0) {$X$};

\node[dot] (x1) at (-0.7,0.2) {};
\node[left] at (x1) {$x$};
\node[dot] (x2) at (0.7,-0.5) {};
\node[dot] (x3) at (0.7,0.4) {};

\draw[arr,blue] (x3) -- (x1);
\draw[arr,blue] (x1) -- (x2);
\draw[arr,blue] (x2) -- (x3);
\node[blue] at (-0.1,0.6) {$\gamma_{xx}$};

\node[dot] (transit1) at (-0.7,-1.2) {};
\draw[arr,red] (x1) -- (transit1);

\draw[cls] (0,-2.4) ellipse (1.6 and 0.9);
\node at (2.5,-2.4) {$Y$};

\node[dot] (y) at (-0.7,-2.4) {};
\node[left] at (y) {$y$};
\draw[arr,red] (transit1) -- (y);
\node[red] at (-1.1,-1.2) {$\gamma'_{xy}$};

\end{tikzpicture}
\end{center}

We set $\gamma_k = (\gamma_{xx})^k \gamma'_{xy} : x \rightsquigarrow y$ and apply Lemma~\ref{infinite_paths}. We obtain that $C$ contains infinitely many intervals.
\end{proof}

\subsection{Local structure of black-and-white classes}

\begin{definition}
Let the edge $v \rightarrowtail w$ of the graph $G'$ correspond to the transition from some half-interval $\nabla_{\alpha_1\ldots\alpha_n}^s$ to $\nabla_{\alpha_1\ldots\alpha_n d}^s$,
($d \in \{0, 1, \ldots, s-1\}$). The number $d$ will be called the \emph{index} of the edge, denoted $\mathrm{ind}(v, w)$.
\end{definition}

Let $v$ be a vertex of some cyclic partial
equivalence class $X$, and let $k$ be the index of the unique edge of the transition within the class to $v' \in X$.
For $v$ we split the set of other immediate descendants into two sets
\[
S^{-}(v) := \{\, w \mid v \rightarrowtail  w,\ \mathrm{ind}(v,w) < k\},\]
\[
S^{+}(v) := \{\, w \mid v \rightarrowtail  w,\ \mathrm{ind}(v,w) > k\}.\]

% add a figure here.

\begin{lemma}\label{Each_S_is_1colored}
For the finiteness of the number of intervals in $C$ it is necessary that for each vertex $v \in K$, where $K$ is a black-and-white class from $G'$, each of the sets $S^{-}(v)$ and $S^{+}(v)$ be monochromatic, that is, consist only of black or only of white vertices.
\end{lemma}
\begin{proof}
By Lemma \ref{necessary_condition1_branching} we assume that all black-and-white classes are cyclic.

Suppose that in some black-and-white class $K$ there exists a vertex $a$ for which the set $S^{-}(a)$ or $S^{+}(a)$ is not monochromatic.
Let this be the set $S^{-}(a)$ (the case $S^+(a)$ is symmetric). There are no black-and-white descendants (Lemma \ref{necessary_condition2_connecting}), so there exist edges $B \colon a \rightarrowtail b$, $W \colon a \rightarrowtail w$, where $b$ is black and $w$ is white, $b, w \in S^{-}(a)$.

Inside the class $K$ we fix some path $\gamma_{aa} \colon a \rightsquigarrow a$. Consider the paths
\[\gamma_i^B \coloneqq (\gamma_{aa})^i B \colon a \rightsquigarrow b; \quad
\gamma_i^W \coloneqq (\gamma_{aa})^i W \colon a \rightsquigarrow w, \qquad i \in \mathbb{N}.\]
No path is a prefix of another: when the shorter path leaves the class $K$ at its last step, the longer one still remains in this class; when the lengths are equal, the paths are different.

For some $\nabla_a$ (with prefix $a$) consider the subcylinders obtained as a result of the transitions $\gamma_i^B$ and $\gamma_i^W$.
%By assumption the set $S^{-}(a)$ is not monochromatic.
The cylinders $\gamma_m^B(\nabla_a)$ and $\gamma_m^W(\nabla_a)$ lie to the left of
$(\gamma_{aa})^{m+1}(\nabla_a)$.
If $\mathrm{ind}(a, b) < \mathrm{ind}(a, w)$, then the disjoint cylinders from left to right have the order
\[ \gamma_1^B(\nabla_a),
\gamma_1^W(\nabla_a),
\gamma_{2}^B(\nabla_a),
\gamma_{2}^W(\nabla_a), \ldots,
\gamma_i^B(\nabla_a),
\gamma_i^W(\nabla_a),\]
therefore $\nabla_a$ contains at least $i$ (that is, infinitely many) intervals from $C$ (one interval in each of the $i$ cylinders with black prefixes, Lemma~\ref{lemma_black}), separated by gaps (cylinders with white prefixes, Lemma~\ref{lemma_gap}).
%Since $i$ is arbitrary, $\nabla_a$ contains infinitely many intervals from $C$.

These arguments remain valid also for $\mathrm{ind}(a, b) > \mathrm{ind}(a, w)$ (only $\gamma_i^B(\nabla_a)$ and $\gamma_i^W(\nabla_a)$ are interchanged).
\end{proof}

\begin{definition}
We say that a vertex $v$ has the \emph{property}:
\begin{itemize}
\item
\emph{$D_1$}, if all vertices from $S^{-}(v)$ are black and all vertices from $S^{+}(v)$ are white;

\item
\emph{$D_2$}, if all vertices from $S^{-}(v)$ are white and all vertices from $S^{+}(v)$ are black.

\end{itemize}
\end{definition}

\begin{remark}\label{empty_s1_or_s2}
If one of the sets $S^{-}(v), S^{+}(v)$ is monochromatic and the other is empty, then one of the properties $D_1$, $D_2$ holds.
\end{remark}

\begin{lemma}\label{all_v_are_D1_or_D2_type}
For the finiteness of the number of intervals in $C$ it is necessary that each vertex from a black-and-white class have property $D_1$ or $D_2$.
\end{lemma}
\begin{proof}
We assume that all black-and-white classes are cyclic (by Lemma~\ref{necessary_condition1_branching}).

Suppose that in some black-and-white class $K$ there exists a vertex $a$ for which neither $D_1$ nor $D_2$ holds. By Lemma \ref{Each_S_is_1colored} $S^{-}(a)$ and $S^{+}(a)$ are monochromatic, so the violation of conditions $D_1$ and $D_2$ means that the colors of $S^{-}(a)$ and $S^{+}(a)$ coincide.
Without loss of generality we assume the color is black. Since the class $K$ is black-and-white and, by Lemma \ref{necessary_condition2_connecting}, has no black-and-white descendants, an edge from it goes to a white vertex. Denote this edge by $W \colon a_2 \rightarrowtail w$, where $a_2 \in K$ and $w$ is white, $w \notin K$. Note that $a_2 \neq a$, since $a$ has only black descendants. We fix an edge
$B \colon a \rightarrowtail b$, where $b$ is a black vertex outside $K$, and we choose it so that $b \in S^+ (a)$ if $w \in S^+ (a_2)$ and $b \in S^- (a)$ if $w \in S^- (a_2)$. Such a choice is possible because $S^+(a) \neq \varnothing$ and $S^-(a) \neq \varnothing$, otherwise by Remark \ref{empty_s1_or_s2} one of the properties $D_1, D_2$ would hold.

Inside the class $K$ we fix a cycle from vertex $a$ that passes through vertex $a_2$, that is, $\gamma_{aa} := \gamma_{a a_2}\,\gamma_{a_2 a}$. Consider the paths
\[ \gamma_i^B \coloneqq (\gamma_{aa})^i\, B \colon a \rightsquigarrow b, \qquad
  \gamma_i^W \coloneqq (\gamma_{aa})^i\,\gamma_{a a_2}\, W \colon a \rightsquigarrow w, \qquad i \in \mathbb{N}. \]

The rest of the proof follows the same scheme as in the previous lemma. No path is a prefix of another.
For some cylindrical interval $\nabla_a$ consider the subcylinders obtained as a result of the transitions $\gamma_i^B$ and $\gamma_i^W$. If $w \in S^- (a_2)$ and, accordingly, $b \in S^- (a)$, then the disjoint subcylinders are arranged in the order
$ \gamma_1^B(\nabla_a),
\gamma_1^W(\nabla_a),
\gamma_{2}^B(\nabla_a),
\gamma_{2}^W(\nabla_a), \ldots.$
If $w \in S^+ (a_2)$ and, accordingly, $b \in S^+ (a)$, then we have the order
$ \gamma_i^W(\nabla_a),
\gamma_i^B(\nabla_a),
\gamma_{i-1}^W(\nabla_a),
\gamma_{i-1}^B(\nabla_a), \ldots$.

Hence $\nabla_a$ contains infinitely many intervals from $C$.
\end{proof}

\begin{lemma}\label{necessary_condition3_middle}
For the finiteness of the number of intervals in $C$ it is necessary that in each black-and-white partial equivalence class all vertices have the same property: either all have property $D_1$, or all have $D_2$.
\end{lemma}
\begin{proof}
We assume that the necessary conditions for the finiteness of the number of intervals (Lemmas~\ref{necessary_condition1_branching},~\ref{necessary_condition2_connecting}~\ref{all_v_are_D1_or_D2_type}) hold: all black-and-white classes are cyclic; a black-and-white class has no black-and-white descendant; all vertices of black-and-white classes have one of the properties $D_1$ or $D_2$.

Suppose that in some black-and-white class $K$ there exist both vertices with property $D_1$ and vertices with property $D_2$. We show that then one can choose two vertices $a_1$ with property $D_1$ and $a_2$ with property $D_2$ such that at least one of the following conditions holds:

\begin{itemize}
\item[(1)]
$S^-(a_1) \neq \varnothing$ and $S^-(a_2) \neq \varnothing$;

\item[(2)]
$S^+(a_1) \neq \varnothing$ and $S^+(a_2) \neq \varnothing$.
\end{itemize}

We split into cases depending on whether there exists in $K$ a vertex with property $D_1$ with $S^+ \neq \varnothing$ and whether there exists a vertex with property $D_1$ with $S^- \neq \varnothing$:

\begin{description}
\item[(A) (yes, yes).] There exists $a_1^+: S^+(a_1^+) \neq \varnothing$ with property $D_1$ and there exists $a_1^-: S^-(a_1^-) \neq \varnothing$ with property $D_1$. Let $a_2$ be some vertex with property $D_2$. If $S^-(a_2) \neq \varnothing$, then the desired pair is $(a_1^-, a_2)$, for which condition (1) holds. If $S^-(a_2) = \varnothing$, then, since $|S^-(a_2) \cup S^+(a_2)| = s - 1 > 0$, we have $S^+(a_2) \neq \varnothing$. The desired pair will be $(a_1^+, a_2)$, for which condition (2) holds.

\item[(B) (yes, no).] There exists $a_1^+: S^+(a_1^+) \neq \varnothing$ with property $D_1$ and there does not exist a vertex $a_1^-: S^-(a_1^-) \neq \varnothing$ with property $D_1$. The absence of vertices of the form $a_1^-$ means that there is no vertex with property $D_1$ with a black outgoing edge. The class $K$ is black-and-white, so there exists a vertex $a_2$ with property $D_2$ that has a black outgoing edge, that is, $S^+(a_2)\neq \varnothing$. The desired pair is $(a_1^+, a_2)$, and condition (2) holds.

\item[(C) (no, yes).] Analogous to (B).

\item[(D) (no, no).] Then for any vertex $a_1$ with property $D_1$ we have $S^-(a_1) = S^+(a_1) = \varnothing$, but all vertices have $s$ outgoing edges, where $s \geq 2$, therefore $S^-(a_1) \cup S^+(a_1) \neq \varnothing$. We have a contradiction, hence this case is impossible.
\end{description}

Thus it is proved that there always exists a pair of vertices $(a_1,a_2)$ such that (1) or (2) holds. Let this be (1), that is, $S^-(a_1) \neq \varnothing$ and $S^-(a_2) \neq \varnothing$.

For the vertex $a_1$ with property $D_1$ we denote the edge with index $0$ by $B \colon a_1 \rightarrowtail b$. Since $S^-(a_1) \neq \varnothing$, we have $b \in S^-(a_1)$. By property $D_1$ all vertices from $S^-(a_1)$ are black, so $b$ is a black vertex outside $K$. Similarly for the vertex $a_2$ with property $D_2$: we denote the edge with index $0$ by $W \colon a_2 \rightarrowtail w$. Then $w \in S^-(a_2) \neq \varnothing$, so $w$ is a white vertex outside $K$.

The class $K$ is strongly connected; inside the class $K$ we fix paths $\gamma_{a_1 a_2}$, $\gamma_{a_2 a_1}$ and denote the cycle
$\gamma_{a_1 a_1} \coloneqq \gamma_{a_1 a_2} \gamma_{a_2 a_1}$.
Consider the paths:
\[ \gamma_i^B \coloneqq (\gamma_{a_1 a_1})^i\, B \colon a_1 \rightsquigarrow b, \qquad
  \gamma_i^W \coloneqq (\gamma_{a_1 a_1})^i\,\gamma_{a_1 a_2}\, W \colon a_1 \rightsquigarrow w, \qquad i \in \mathbb{N}. \]

For some $\nabla_{a_1}$ consider the subcylinders obtained as a result of the transitions $\gamma_i^B$ and $\gamma_i^W$. As in the previous lemmas, no path is a prefix of another, so these subcylinders are disjoint.

Inside the subcylinder $(\gamma_{a_1 a_1})^i (\nabla_{a_1})$ the edge $B$ has index $0$ and leads to the leftmost subcylinder, while the first step of $\gamma_{a_1 a_2}$ has index greater than $0$, so it leads to a subcylinder that lies to the right of $\gamma_i^B(\nabla_{a_1})$.
Similarly, inside the subcylinder $\bigl((\gamma_{a_1 a_1})^i\gamma_{a_1 a_2}\bigr)(\nabla_{a_1})$ the edge $W$ has index $0$ and leads to the leftmost subcylinder, while the first step of $\gamma_{a_2 a_1}$ leads to a subcylinder that lies to the right of $\gamma_i^W(\nabla_{a_1})$. Therefore the disjoint subcylinders, when considered from left to right, are arranged in the order
$ \gamma_1^B(\nabla_{a_1}),
\gamma_1^W(\nabla_{a_1}),
\gamma_{2}^B(\nabla_{a_1}),
\gamma_{2}^W(\nabla_{a_1}), \ldots.$
Since $i$ is arbitrary, $\nabla_{a_1}$ contains infinitely many intervals.

The case of condition (2) is symmetric to the case of condition (1), therefore the arguments in case (2) are analogous.
\end{proof}

\begin{lemma}\label{lemma_one_segment}
Let all vertices of a cyclic black-and-white class $K$ have the same property $D_1$ or $D_2$.
Then each cylindrical half-interval with prefix $q \in K$ contains exactly one (maximal within the cylinder) interval/half-interval from $C$.
\end{lemma}

\begin{proof}
We carry out the proof for the case $D_1$. The case $D_2$ is symmetric and is proved analogously.

For an arbitrary vertex $v \in K$ denote by $k = k(v)$ the index of the unique internal edge from it, let $v \rightarrowtail v'$, $v' \in K$. Property $D_1$ means that the vertices in $S^-(v)$ are black and those in $S^+(v)$ are white. Therefore the cylinder $\nabla_{v}$ from left to right is partitioned into three consecutive parts:
\begin{itemize}
  \item subcylinders from $S^-(v)$ (black prefixes) --- consecutive half-intervals that lie entirely in $C$ (Remark~\ref{remark_black_is_full}). Together they form a half-interval that lies in $C$;
  \item the subcylinder $\nabla_{v'}$, which also has a prefix from $K$;
  \item subcylinders from $S^+(v)$ (white prefixes) --- consecutive half-intervals, each of which contains no intervals from $C$ (Remark~\ref{remark_white_no_segment}). Together they form a half-interval that contains no intervals from $C$.
\end{itemize}

% TikZ: a solid black bar on the left, an arrow ``deeper into $C$'', empty in the middle, a wavy bar on the right.

Put $u_0 \coloneqq q$ and fix a cylindrical interval $\nabla_{u_0} = [L_0, R_0)$. Denote the sequence of vertices obtained by iterative transitions along the internal edge by $u_1, u_2, \ldots$, that is, $u_i \overset{k(u_i)}\rightarrowtail u_{i+1}, i \geq 0$.

To the sequence $u_0, u_1, \ldots ,u_i, \ldots$ from $C$ there corresponds a sequence of nested cylinders $\nabla_{u_i} = [L_i, R_i)$, for which, as shown above,
$[L_i, L_{i+1}) \subset C$, while $[ R_{i+1}, R_i )$ contains no interval from $C$.

The sequence $(L_i)$ is nondecreasing, $(R_i)$ is nonincreasing, $R_i - L_i \to 0 (i \to \infty)$.
Hence they tend to some common point $p \in [L_0, R_0]$. Combining the left and right parts, we have $[L_0, p) \subset C$, while $(p, R_0)$ contains no interval from $C$.
We have $p \neq L_0$ and $p \neq R_0$, otherwise the class $K$ would not be black-and-white.
%Suppose that $L_0 = p$. Then $L_i = L_0$ for all $i$, all indices of internal edges equal $0$, therefore all external edges for each of the vertices $(u_i)$ lie in $S^+(u_i)$. The class $K$ is cyclic, therefore the sequence $(u_i)$ includes all vertices of the class. Thus by $D_1$ all external edges from the class $K$ lead to white vertices. Hence there is no path from the class $K$ to a black vertex, which contradicts the black-and-white character of $K$. Analogously, the assumption $p = R_0$ would exclude a path to a white vertex. Hence $L_0 < p < R_0$.
Hence $\nabla_q$ contains exactly one maximal (within the cylinder) interval $[L_0, p]$ from $C$ ($p \in C$, since $C$ is closed).
\end{proof}

\begin{theorem}[A criterion for the finiteness of the number of intervals in $C$]\label{theorem_finiteness_criterion}
The set $C$ contains only finitely many intervals if and only if the following conditions hold simultaneously:

1) all black-and-white classes are cyclic,

2) there are no paths between black-and-white classes,

3) for a black-and-white class all vertices in it have the same property: either all have property $D_1$, or all have $D_2$.
\end{theorem}

\begin{proof}
\emph{Necessity.} Proved in the three preceding lemmas (Lemma~\ref{necessary_condition1_branching}, Lemma~\ref{necessary_condition2_connecting}, Lemma~\ref{necessary_condition3_middle}).

\emph{Sufficiency.} Let conditions (1--3) hold.

We shall say that a path $a \rightsquigarrow b$ is a path to some class $K \ni b $ if $b$ is its first (unique) entry into the class $K$.

Suppose that a path $p \rightsquigarrow K_1$ from some starting vertex $p$ to a black-and-white class $K_1$ passes through some other class $K_0$, that is, $p \rightsquigarrow K_1 = p \rightsquigarrow K_0 \rightsquigarrow K_1$. Since $K_1$ is black-and-white, by Lemma~\ref{lemma_color_inheritance} the class $K_0$ is also black-and-white. This contradicts condition (2). Hence paths from starting vertices to black-and-white classes pass only through transient vertices.
Transient vertices do not belong to any class, and hence do not belong to any cycle. The graph $G'$ has finitely many transient vertices (in general there are finitely many vertices). There are also finitely many black-and-white classes and starting vertices, therefore there exist only finitely many paths from starting vertices to black-and-white classes.

Each path $p \rightsquigarrow K$, where $p$ is a starting vertex and $K$ is a black-and-white class, determines a certain cylindrical half-interval $\nabla_X$, which by Lemma~\ref{lemma_one_segment} contains exactly one interval from $C$ within its bounds.
Analogously, the number of paths $p \rightsquigarrow B$, where $p$ is a starting vertex and $B$ is a black class, that pass only through transient vertices in the graph $G'$, is also finite. Each path $p \rightsquigarrow B$ determines a certain cylindrical half-interval $\nabla_Y \subset C$ (Remark~\ref{remark_black_is_full}).

Two families of cylindrical half-intervals have been described, each of which contains exactly one (half-)interval from $C$ within its bounds. There are finitely many of them, therefore the union of all these cylinders, which we denote by $F$, contains only finitely many intervals from $C$.

Suppose there exists an interval $I \subset C$ such that $I \cap F = \varnothing$. Inside $I$ we choose some cylinder $\nabla_u$; $\nabla_u \subset I \subset C$, therefore the prefix $u$ is black (Remark~\ref{remark_black_is_full}).
Consider the path $p \rightsquigarrow u$.

Suppose the path $p \rightsquigarrow u$ does not consist only of transient vertices. Denote by $u_1$ the first vertex on it that belongs to some partial equivalence class. It has a descendant $u$, therefore by Lemma~\ref{lemma_color_inheritance} it is black or black-and-white, and hence falls into one of the families described above.
The cylindrical half-interval corresponding to the path $p \rightsquigarrow u_1$ contains $\nabla_u$, a contradiction with the assumption.

Suppose the path $p \rightsquigarrow u$ consists only of transient vertices. Since the prefix $u$ is black, there exists a path $u \rightsquigarrow t$, where $t$ is a terminal black class. Denote by $u_1$ the first vertex on it that belongs to some partial equivalence class. Since $u$ is black, all of its descendants are black, hence $u_1$ is black.
The path $p \rightsquigarrow u\rightsquigarrow u_1$ determines a cylindrical interval that is an element of the second family, a contradiction.

Thus the two families described exhaust all intervals from $C$.

Hence $C$ contains only finitely many intervals.
\end{proof}

\section{Examples of applications}

For multigeometric series $(k_1, k_2, \ldots, k_m; \frac{1}{s})$ the criterion for the existence of an interval and the criterion for the finiteness of their number in the set of subsums make it possible to use the Guthrie--Nymann classification theorem and to determine the topological structure of the set of subsums; that is, for multigeometric series of the form~\eqref{multigeom_series} a criterion for being a Cantorval has been obtained.

Not every set $C$ is the set of incomplete sums of a series,
therefore the Guthrie--Nymann theorem is not directly applicable to it. %, that is, the existence of infinitely many intervals does not imply that $C$ is a Cantorval.

\emph{Example 1}. The Guthrie--Nymann series: $A = \{0,2,3,5\}, s=4$. The graph $G$:

{\centering
\begin{tikzpicture}[
  vertex/.style={circle, draw, minimum size=10pt},
  arr/.style={->, line width=1pt}
]

% ---- vertices
\node[vertex] (v00) at (-4,0) {$(0,0)$};
\node[vertex] (v10) at (-1.5,0) {$(1,0)$};
\node[vertex] (v01) at (1.5,0) {$(0,1)$};
\node[vertex] (v11) at (4,0) {$(1,1)$};

% (0,0)
\draw[arr,->] (v00) to[out=60, in=120, looseness=3]   (v00);
\draw[arr,->] (v00) to[out=55, in=125, looseness=4]   (v00);
\draw[arr,->] (v00) to[out=50, in=130, looseness=5]   (v00);
\draw[arr,->] (v00) to[out=45, in=135, looseness=6]   (v00);

% (0,1)
\draw[arr,->] (v01) to[out=60, in=120, looseness=3]   (v01);
\draw[arr,->] (v01) to[out=45, in=135, looseness=4]   (v01);
\draw[arr, bend left=15] (v01) to (v10);
\draw[arr, bend left=15] (v01) to (v11);

% (1,0)
\draw[arr,->] (v10) to[out=60, in=120, looseness=3]   (v10);
\draw[arr,->] (v10) to[out=45, in=135, looseness=4]   (v10);
\draw[arr, bend left=15] (v10) to (v01);
\draw[arr, bend left=15] (v10) to (v00);

% (1,1)
\draw[arr,->] (v11) to[out=60, in=120, looseness=3]   (v11);
\draw[arr,->] (v11) to[out=55, in=125, looseness=4]   (v11);
\draw[arr,->] (v11) to[out=50, in=130, looseness=5]   (v11);
\draw[arr,->] (v11) to[out=45, in=135, looseness=6]   (v11);

\end{tikzpicture}
\par}

White class $(0,0)$, one black class $(1,1)$.
Starting vertices: $(0,1)$ for the interval $[0,1)$ and $(1,0)$ for $[1,2)$; they form a black-and-white class $\{ (0,1), (1,0)\}$. There exists a path from a starting vertex to a black class, hence the set of subsums contains an interval. The black-and-white class is not cyclic, hence there are infinitely many intervals. By the Guthrie--Nymann classification theorem the set of subsums of this series is a Cantorval.

\emph{Example 2}. We present a new multigeometric series for which the Cantorval property of the set of subsums is not verified by the already known sufficient conditions described in the papers~\cite{Bartoszewicz_Filipczak_Szymonik_2014} and~\cite{Nowakowski_OnANewCondition}. If one restricts oneself to series of the form (\ref{multigeom_series}), then these sufficient conditions, when passing to the graph, mean reaching the terminal black class $(1,1,\ldots,1)$. For the series given below the vertex $(1,1,\ldots,1)$ is not reachable from the starting prefixes.

The series $(1,4,16; \tfrac{1}{5})$: $V = \{0,1,4,5,16,17,20,21\}$, $M = 6$. The condensed graph $\bar{G'}$:

{\centering
\begin{tikzpicture}[
  scale=1.0,
  vertex/.style={ellipse, draw, inner sep=2pt, font=\scriptsize},
  cls/.style={draw, dashed, line width=0.9pt, rounded corners=9pt, inner sep=5pt, align=center, font=\scriptsize},
  arr/.style={->, line width=1pt}
]

% ---- starting transient vertices
\node[vertex] (T1) at (-5.1,4.3) {\texttt{000100}};
\node[vertex] (T2) at (-2.7,4.3) {\texttt{001000}};
\node[vertex] (T3) at (0.9,4.3)  {\texttt{000010}};
\node[vertex] (T4) at (3.7,4.3)  {\texttt{010000}};

\node at (-5.1,4.8) {\scriptsize $[2,3)$};
\node at (-2.7,4.8) {\scriptsize $[3,4)$};
\node at (0.9,4.8)  {\scriptsize $[1,2)$};
\node at (3.7,4.8)  {\scriptsize $[4,5)$};

% black-and-white class
\node[cls] (C1) at (-0.7,3.0) {\itshape black-and-white class\\[1pt]
\texttt{000001}\ \texttt{001100}\ \texttt{000110}\\
\texttt{100110}\ \texttt{110000}\ \texttt{100000}};

% transient vertices
\node[vertex] (t5) at (-3.4,2.0) {\texttt{000011}};
\node[vertex] (t6) at (2.6,2.0)  {\texttt{011000}};

% black class
\node[cls] (K1) at (5,1.5) {\itshape black class\\[1pt] \texttt{011001}};

% terminal classes
\node[cls] (W) at (-4.8,-0.8) {\itshape white terminal class\\[1pt] \texttt{000000}};

\node[cls] (B) at (1.0,-0.7) {\itshape black terminal class\\[1pt]
\texttt{110011}\ \texttt{100111}\ \texttt{001110}\ \texttt{011100}\\
\texttt{111001}\ \texttt{011101}\ \texttt{111011}\ \texttt{110111}\\
\texttt{101110}\ \texttt{001111}\ \texttt{011110}\ \texttt{111100}\\
\texttt{101111}\ \texttt{111101}};

% edges
% 000100
\draw[arr] (T1) to (C1);
\draw[arr] (T1) to (W);
% 001000
\draw[arr] (T2) to (C1);
\draw[arr] (T2) to (t5);
\draw[arr] (T2) to[bend right=25] (W);
% 000010
\draw[arr] (T3) to (C1);
\draw[arr] (T3) to (t6);
\draw[arr] (T3) to[bend left=20] (B);
% 010000
\draw[arr] (T4) to (C1);
\draw[arr] (T4) to (t6);
\draw[arr] (T4) to (K1);
\draw[arr] (T4) to[bend left=20] (B);
% black-and-white class
\draw[arr] (C1) to (t5);
\draw[arr] (C1) to (t6);
\draw[arr] (C1) to (W);
\draw[arr] (C1) to[bend right=15] (K1);
\draw[arr] (C1) to (B);
% transient vertices
\draw[arr] (t5) to (B);
\draw[arr] (t6) to (K1);
\draw[arr] (t6) to (B);
% black class
\draw[arr] (K1) to (B);

\end{tikzpicture}
\par}

From the starting vertices a black terminal class is reachable. At the same time the black-and-white class is not cyclic: for example, from the vertex $001100$ there emanate four internal edges.
Hence the set of subsums of the series contains infinitely many intervals, that is, it is a Cantorval.

\section*{Acknowledgements}
The English translation from Ukrainian original text was done with usage of an AI language model Claude Fable 5.1. We have checked and edited the translation and take full responsibility for the content.

\bibliographystyle{plain}
\bibliography{mybiblio-ed2}
\end{document}

%%%%%%%%%%%%%%%% Insert metadata of the paper

\selectlanguage{ukrainian}
%\selectlanguage{english}

% Set the title of the paper:
% If yout title is long and will not fit to hte top margin of the paper
% then specify a short title in square brackets
\title[Канторвальність при натуральній основі]{Критерій канторвальності множин підсум одного класу мультигеометричних рядів}

\author[М. Працьовитий]{Микола Працьовитий}
\organization{Інститут математики НАН України, УДУ імені Михайла Драгоманова}
\email{prats4444@gmail.com}
\orcid{0000-0001-6130-9413}
%\thanks{This work was supported by a grant from the Simons Foundation (SFI-PD-Ukraine-00014586, M.P., O.V.).}
%%%%%%%%%%%% information on the second author:

\author[Я. Виннишин]{Ярослав Виннишин}
\organization{Київський природничо-науковий ліцей №~145}
\email{vynnyar@gmail.com}

\author[О. Виннишин]{Олег Виннишин}
\organization{Інститут математики НАН України}
\email{oleh.vynnyshyn@imath.kiev.ua}
\orcid{0009-0006-7113-1501}

\abstract{ukrainian}{У роботі досліджується топологічна структура множини підсум (неповних сум) мультигеометричних рядів вигляду $\frac{k_1}{s} + \frac{k_2}{s} + \ldots + \frac{k_m}{s} + \frac{k_1}{s^2} + \frac{k_2}{s^2} + \ldots+ \frac{k_m}{s^2} + \frac{k_1}{s^3} + \frac{k_2}{s^3} +\ldots$, де $k_1, \ldots, k_m, s \in \mathbb{N}$.
У термінах теорії графів встановлено необхідні та достатні умови того, що множина підсум даного ряду не є скінченним об'єднанням відрізків, а також необхідні та достатні умови того, що ця множина містить відрізок. Наслідком цих тверджень є критерій канторвальності множини підсум ряду.
Наведено приклад застосування критерію до ряду, канторвальність якого раніше не була відомою.

 Більш загальним об'єктом вивчення топологічних властивостей є множини чисел, на зображення елементів яких в системі з натуральною основою та надлишковим алфавітом накладені заборони на використання цифр. Саме для цього об'єкта проведено обґрунтування вказаних тверджень.
}

% \dedication{...}
\thanks{This work was supported by a grant from the Simons Foundation (SFI-PD-Ukraine-00014586, P.M., V.O.).}
\keywords{Канторвал; множина неповних сум; мультигеометричний ряд; критерій канторвальності; множина канторівського типу; система числення з натуральною основою та надлишковим алфавітом
}
\msc{40A05, 28A80, 11A67, 05C90}
%
% \doi{}
% \received{}
% \accepted{}

\maketitle

%%% Text of the paper

\section{Вступ}
Нехай $r$ --- фіксоване натуральне число, $A_r = \{0,1, \ldots, r\}$ --- алфавіт; $L_r = A_r \times A_r \times A_r \times \cdots$ --- множина послідовностей елементів алфавіту, зокрема $L_1$ --- множина послідовностей з нулів та одиниць; $x_1 + x_2 + \cdots + x_k+ \cdots$ --- абсолютно збіжний числовий ряд з дійсними відмінними від нуля членами. 

Нагадаємо~\cite{Kakeya_1914_partial_sums}, що множиною підсум (неповних сум) даного ряду називається множина

\[ \mathcal{E}\bigl(x_n\bigr) \coloneqq 
\left\{ x \; : \;
x = \sum_{n=1}^{\infty} c_n x_n, \quad
(c_n) \in L_1 \right\}. \]

Вперше топологічні властивості множин неповних сум абсолютно збіжних рядів досліджував Соічі Какея~\cite{Kakeya_1914_partial_sums}, який довів, що вони є досконалими; знайшов необхідні та достатні умови, за яких множина підсум є відрізком або скінченним об'єднанням відрізків, а також знайшов достатні умови ніде не щільності множини підсум. Пізніше, у роботі~\cite{3_Guthrie_Nymann} було отримано класифікаційний результат:

\begin{theorem}[Класифікаційна теорема Гатрі-Німана]
\label{Hatri_Nymann_th}
Множина неповних сум (підсум) абсолютно збіжного ряду є

\begin{enumerate}
%\item 
%або скінченною множиною
\item 
або об'єднанням скінченної кількості відрізків
\item
або гомеоморфною множині Кантора
\item 
або канторвалом, тобто множиною, гомеоморфною множині 
\[C \,\cup\, \bigcup_{n=1}^{\infty} G_{2n-1},\]
де $C$ --- класична канторова множина, а $G_n$ --- відрізки,
що видаляються на $n$-му кроці її побудови.
\end{enumerate}
\end{theorem}

Поняття канторвалу було введено в роботі~\cite{Mendes_Oliveira_MCantorvals_1994}, найбільш відомим прикладом є канторвал Гатрі-Німана~\cite{3_Guthrie_Nymann}, що є множиною неповних сум ряду з членами $x_{2n-1} = \frac{2}{4^n}, x_{2n} = \frac{3}{4^n}$. 

Мультигеометричні ряди вигляду
%\[ \sum_{n=1}^{\infty} \sum_{j=1}^{m} k_j q^n, \]
\[k_1 q + k_2 q + \ldots + k_m q + k_1 q^2 + k_2 q^2 + \ldots+ k_m q^2 + k_1 q^3 + \ldots,\]
де $ k_1,\ldots,k_m > 0, \quad q \in (0,1)$
вивчались в роботах
~\cite{Ferens1984,
Nitecki2015, 
Bartoszewicz_Filipczak_Szymonik_2014, 
Karvatskyi_Murillo_Viruel_2024_generalized_multigeom, 
Pratsiovytyi_Karvatskyi_2023_trigonometric, 
PratsiovytyiKarvatskyiMakarchuk2026} 
та інших. Для них у роботі~\cite{Bartoszewicz_Filipczak_Szymonik_2014} наведено достатні умови канторвальності множини неповних сум,  введено позначення мультигеометричного ряду як $(k_1, k_2, \ldots, k_m; q)$, яке ми далі будемо використовувати.
В роботі~\cite{Nowakowski_OnANewCondition} отримано нову достатню умову канторвальності множини підсум.

Сучасний стан і перелік відкритих проблем стосовно тополого-метричних властивостей множин неповних сум ряду збіжних рядів наведено в оглядовій праці~\cite{GlabPrusWisniowski2025}, де однією з центральних проблем вказано критерій канторвальності множини неповних сум ряду. 

Об'єктом дослідження у даній роботі є 
ряд 

\begin{equation}\label{multigeom_series}
(k_1, k_2, \ldots, k_m; \frac{1}{s}) = \frac{k_1}{s} + \ldots + \frac{k_m}{s} + \frac{k_1}{s^2} + \ldots+ \frac{k_m}{s^2} + \ldots,
\end{equation}
де $s,k_1, \ldots, k_m$ --- фіксовані натуральні числа, $s \geq 2$, а також множини з ним пов'язані. 
Історично саме мультигеометричні ряди з цілими коефіцієнтами $k_1, \ldots, k_m$ дали приклади множин підсум, які не є ні скінченними об'єднаннями відрізків, ні множинами канторівського типу~\cite{3_Guthrie_Nymann,Ferens1984, Jones2011}. 
Також множини підсум мультигеометричних рядів з цілими коефіцієнтами природно виникають і в теорії ймовірностей ~\cite{PratsiovytyiKarvatskyiMakarchuk2026, PratsiovytyiRatushniak2025}: вони є спектрами розподілів відповідних нескінченних згорток Бернуллі. 

Отриманий у цій статті критерій канторвальності для щойно описаного класу мультигеометричних рядів (\ref{multigeom_series}) характеризує топологічний тип множин їх підсум, тож його можна розглядати як часткове вирішення задачі про топологічну структуру множини підсум ряду.

\section{Множини чисел з обмеженнями на використання цифр у своїх зображеннях у системі з натуральною основою та надлишковим алфавітом}

%\subsection{Системи числення з надлишковим алфавітом}

Далі $s$ та $r$ --- фіксовані натуральні числа, $2 \leq s \leq r$. Відомо~\cite{Renyi1957}, що для будь-якого числа $x \in [0, \frac{r}{s-1}]$ існує принаймні одна послідовність $(\alpha_k) \in L_r$, така що $x = \sum_{k=1}^{\infty}\alpha_ks^{-k } = \Delta_{\alpha_1 \alpha_2 \ldots \alpha_n \ldots}^{r_s}$.
Останній символічний запис $\Delta_{\alpha_1 \alpha_2 \ldots \alpha_n \ldots}^{r_s}$ 
ми називаємо $r_s$-зображенням числа $x$, $\alpha_n$ --- $n$-тою цифрою цього зображення. Оскільки алфавіт $A_r$ надлишковий, то, взагалі кажучи, числа відрізка $[0, \frac{r}{s-1}]$ можуть мати єдине зображення, скінченну, зліченну і навіть континуальну множину різних зображень~\cite{PratsiovytyiVynnyshyn2026, Pratsiovytyi_Ratushniak_2026}.

Основним об'єктом розгляду у даній роботі є множина 
\[C = C[r_s,V] = \left\{ x : x = \Delta_{\alpha_1 \alpha_2 \ldots \alpha_n \ldots}^{r_s}, \alpha_n \in V \subset A_r \right\},\]
де $V$ --- фіксована власна не одноелементна підмножина алфавіту $A_r$. У даній роботі нас цікавлять топологічні властивості множини $C$ за умови $\{0,r\} \subset V$.

Для мультигеометричного ряду $(k_1, k_2, \ldots, k_m; \frac{1}{s})$, де $k_i, s$ --- натуральні числа, множина неповних сум $\mathcal{E}(x_n)$ є множиною $C[r_s, V]$ за умови $V = \{v : v = \sum_{i=1}^m c_i k_i, \quad (c_1,c_2, \ldots , c_m) \in A_1^m \}$.

\emph{Приклад 1}. Ряд $(3,2; \tfrac{1}{4})$ (Гатрі-Німана~\cite{3_Guthrie_Nymann}). Для нього $V = \{0,2,3,5\}$, а множина неповних сум $C$ --- канторвал.

\emph{Приклад 2}. Ряд $(6,1; \tfrac{1}{4})$ (приклад Кеньйона~\cite{Nitecki2015}). В даному випадку $V = \{0,1,6,7\}$, $C$ ---канторвал.

\emph{Приклад 3}. Ряд $(1,8; \tfrac{1}{4})$~\cite{GlabPrusWisniowski2025}. Алфавіт $V = \{0,1,8,9\}$, $C = [0,1] \cup [2,3]$.

Не всі множини $C[r_s,V]$ є множиною підсум деякого мультигеометричного ряду вигляду \eqref{multigeom_series}. Також, взагалі кажучи, для таких множин класифікаційна теорема Гатрі-Німана не застосовується.

\section{Структура множини $C[r_s,V]$}

Далі використовуватимемо позначення $\nabla_{\alpha_1 \alpha_2 \ldots \alpha_m}^s $ для $s$-кового циліндричного піввідрізка  рангу $m$ з основою  $\alpha_1 \alpha_2 \ldots \alpha_m$, тобто 

\[ \nabla_{\alpha_1 \alpha_2 \ldots \alpha_m}^s := 
\left[ \sum_{k=1}^{m} \frac{\alpha_k}{s^k}, \;
\sum_{k=1}^{m} \frac{\alpha_k}{s^k} + \frac{1}{s^m} \right), \]
де $\alpha_k \in \{0,1,\ldots,s-1\}$.

Очевидно, що 
\[ \nabla_{\alpha_1 \alpha_2 \ldots \alpha_m}^s = 
\nabla_{\alpha_1 \alpha_2 \ldots \alpha_m0}^s \cup
\nabla_{\alpha_1 \alpha_2 \ldots \alpha_m1}^s \cup \ldots 
\cup \nabla_{\alpha_1 \alpha_2 \ldots \alpha_m[s-1]}^s
\]
%для циліндричного напіввідрізка рангу $m$ з основою $\alpha_1 \alpha_2 \ldots \alpha_m$. 

Напіввідкритий інтервал $[0;1)$ назвемо циліндричним піввідрізком нульового рангу. Він відповідно розбивається на $s^n$ попарно неперетинних циліндричних піввідрізків $n$-го рангу довжини $\frac{1}{s^n}$.

Для зручності позначимо
\[ {\Delta}_{\alpha_1 \alpha_2 \ldots \alpha_n} \coloneqq 
C[r_s,V] \cap 
\nabla_{\alpha_1 \alpha_2 \ldots \alpha_n}^s.\]
%Через $A \oplus B \coloneqq \{a + b \mid a \in A,\ b \in B\}$ позначатимемо арифметичну суму множин $A$ та $B$.
Також нехай $[0,M)$ --- мінімальний піввідрізок, що містить множину $C[r_s, V]$, враховуючи що $0 \in V$. Зрозуміло, що $M = [\frac{r}{s-1}] + 1$. 

Якщо ввести позначення множин
\[ E_n \coloneqq \left\{ x \; : \;
x = \sum_{k=1}^{n} \frac{\alpha_k}{s^k}, \quad
(\alpha_1, \ldots, \alpha_n) \in V^n \right\}, \]

\[ R_n \coloneqq \left\{ x \; : \;
x = \sum_{k=n+1}^{\infty} \frac{\alpha_k}{s^k}, \quad
\alpha_k \in V, \forall k \in \mathbb{N} \right\}, \]

то, очевидно, що $E_n \subset E_{n+1} \subset C, \, R_{n+1} \subset R_n \subset C \quad \forall n \in \mathbb{N}$, зокрема 
\[C = E_n \oplus R_n, \quad C = [\bigcup_{k=1}^{\infty} E_k] = [\lim_{k \to \infty} E_k],  \]
де $[B]$ --- замикання множини $B$.

Елементи множини $E_n$ називатимемо \emph{точками рангу $n$}. 
\begin{lemma}\label{points_rank_n}
 Кожну точку рангу $n$ можна подати у вигляді $\frac{a}{s^n}$, де $a$ --- деяке ціле невід'ємне число. 
\end{lemma} 
\begin{proof}
Дійсно, кожен елемент множини $E_n$ має вигляд
\[\frac{\alpha_1}{s} + \frac{\alpha_2}{s^2} + \ldots + \frac{\alpha_n}{s^n} = \frac{\sum_{i=1}^n \alpha_is^{n-i}}{s^n},\]
де $\alpha_i \in V$ --- цілі числа.
\end{proof}

Множиною точок рангу $0$ вважатимемо множину з однієї точки $\{0\}$.

\medskip

Розглянемо асоційовану з циліндричним піввідрізком $\nabla_{\alpha_1 \alpha_2 \ldots \alpha_n}^s$ $M$-точкову множину $T$, яка містить числа, означені рівністю
 
\[t_j := \frac{\alpha_1}{s} + \ldots + \frac{\alpha_n}{s^n} - \frac{M-1}{s^n} + \frac{j}{s^n}, \quad j \in \{0,1,\ldots, M-1\}.\] 

\medskip
\begin{tikzpicture}[
  scale=1.05,
  axis/.style={->, line width=1pt},
  dot/.style={circle, fill, inner sep=1.5pt},
  hollow/.style={circle, draw, inner sep=2pt},
  cyl/.style={line width=2pt, black}
]

% вісь
\draw[axis] (0,0) -- (8,0);

% циліндр
\draw[cyl] (6,0.12) -- (7,0.12);
\draw[dashed] (6,0) -- (6,0.12);
\draw[dashed] (7,0) -- (7,0.12);
\node[above] at (6.5,0.18) {$\nabla_{\alpha_1 \alpha_2 \ldots \alpha_n}^s$};

% префіксні точки
\node[dot] at (0,0) {};
\node[dot] at (1,0) {};
%\node[hollow] at (2,0) {};
\node[dot] at (2,0) {};
\node[dot] at (4,0) {};
%\node[hollow] at (5,0) {};
\node[dot] at (5,0) {};
\node[dot] at (6,0) {};

% підписи
\node[below] at (0,-0.1) {\scriptsize $t_0$};
\node[below] at (1,-0.1) {\scriptsize $t_1$};
\node[below] at (2,-0.1) {\scriptsize $t_2$};

\node[below] at (4,-0.1) {\scriptsize $t_{M-3}$};
\node[below] at (5,-0.1) {\scriptsize $t_{M-2}$};
\node[below] at (6,-0.1) {\scriptsize $t_{M-1}$};

\node at (3,0) {$\cdots$};

% відстань 1/s^n
\draw[
  decorate,
  decoration={brace, amplitude=6pt}
] (1,0.1) -- (2,0.1)
node[above] at (1.5,0.2) {$\frac{1}{s^n}$};

\end{tikzpicture}

Точка $t_{M-1}$ є лівим кінцем циліндричного піввідрізка $\nabla_{\alpha_1 \alpha_2 \ldots \alpha_n}^s$, а точки $t_{j-1}$ отримуються зсувом точок $t_j$ на $\frac{1}{s^n}$ вліво.

\begin{lemma}\label{use_only_points_t_i}
Для довільного $\nabla_{\alpha_1 \alpha_2 \ldots \alpha_n}^s$ має місце рівність 

\[
\Delta_{\alpha_1 \alpha_2 \ldots \alpha_n}
=
\bigl( (E_n \cap T)
\oplus R_n \bigr) 
\cap \nabla_{\alpha_1 \alpha_2 \ldots \alpha_n}^s.
\]
\end{lemma}
\begin{proof}
За означенням 
\[\Delta_{\alpha_1 \alpha_2 \ldots \alpha_n} = 
C \cap \nabla_{\alpha_1 \alpha_2 \ldots \alpha_n}^s = 
\bigl( E_n \oplus R_n \bigr) \cap 
[t_{M-1}, t_{M-1} + \tfrac{1}{s^n}). \]
Зазначимо, що для довільних $A \subset \mathbb{R}$, $B \subseteq [b_1, b_2)$
\[\bigl( A \oplus B\bigr) \cap [x,y) =
\Bigl( \bigl(A \cap (x-b_2,\, y - b_1)\bigr)\oplus B \Bigr)
\cap [x,y). \]
Оскільки множина $R_n \subset[0,\frac{M}{s^n})$, то
\[\Delta_{\alpha_1 \alpha_2 \ldots \alpha_n} = 
\Bigl( 
    \bigl(
        E_n \cap (t_{M-1} - \tfrac{M}{s^n},\, 
        t_{M-1} + \tfrac{1}{s^n})
    \bigr)
    \oplus R_n 
\Bigr) 
\cap [t_{M-1}, t_{M-1} + \tfrac{1}{s^n}).\]
 За означенням $t_0 - \tfrac{1}{s^n} = t_{M-1} - \frac{M}{s^n}$, маємо
\[
\Delta_{\alpha_1 \alpha_2 \ldots \alpha_n} = \bigl( (E_n \cap(t_0 - \tfrac{1}{s^n}, 
t_{M-1} + \tfrac{1}{s^n}))
\oplus R_n \bigr) 
\cap \nabla_{\alpha_1 \alpha_2 \ldots \alpha_n}^s.\]

Кожна з точок $t_0, \ldots, t_{M-1}$ може як належати множині $E_n$, так і не належати їй, при цьому жодна інша точка з $(t_0 - \tfrac{1}{s^n}, t_{M-1} + \tfrac{1}{s^n})$ не належить $E_n$ (лема~\ref{points_rank_n}). Тому
\[\Delta_{\alpha_1 \alpha_2 \ldots \alpha_n} = 
\bigl( (E_n \cap \{t_0, t_1, \ldots, t_{M-1}\})
\oplus R_n \bigr) 
\cap \nabla_{\alpha_1 \alpha_2 \ldots \alpha_n}^s. 
\qedhere \]
\end{proof}

\section{Множини префіксних точок та префікси}

\begin{definition}
\emph{Множиною префіксних точок} циліндричного піввідрізка 
$\nabla_{\alpha_1 \ldots \alpha_n}^s$ називається множина
\[ B(\nabla_{\alpha_1 \ldots \alpha_n}^s) := 
\{t_0, t_1, \ldots, t_{M-1}\} \cap E_n. \]
\end{definition}

\begin{definition}
\emph{Префіксом} циліндричного піввідрізка $\nabla_{\alpha_1\ldots\alpha_n}^s$ рангу $n$ називається впорядкований набір $(p_0, p_1, \ldots, p_{M-1})$, де
$p_i=1$, якщо $t_i\in E_n$, і $p_i=0$, якщо $t_i\notin E_n$.
\end{definition}

%Тобто множина префіксних точок задає префікс:
%\begin{tikzpicture}[
%  scale=1.05,
%  axis/.style={->, line width=1pt},
%  dot/.style={circle, fill, inner sep=1.5pt},
%  hollow/.style={circle, draw, inner sep=2pt},
%  cyl/.style={line width=2pt, black}
%]

%% вісь
%\draw[axis] (0,0) -- (8,0);

%% циліндр
%\draw[cyl] (6,0.12) -- (7,0.12);
%\draw[dashed] (6,0) -- (6,0.12);
%\draw[dashed] (7,0) -- (7,0.12);

% префіксні точки
%\node[dot] at (0,0) {};
%\node[dot] at (1,0) {};
%\node[hollow] at (2,0) {};

%\node[dot] at (4,0) {};
%\node[hollow] at (5,0) {};
%\node[dot] at (6,0) {};

%texts
%\node[below] at (0,-0.1) {\scriptsize $t_0$};
%\node[below] at (1,-0.1) {\scriptsize $t_1$};
%\node[below] at (2,-0.1) {\scriptsize $t_2$};

%\node[below] at (4,-0.1) {\scriptsize $t_{M-3}$};
%\node[below] at (5,-0.1) {\scriptsize $t_{M-2}$};
%\node[below] at (6,-0.1) {\scriptsize $t_{M-1}$};

%\node[above] at (0,0.1) {\scriptsize $1$};
%\node[above] at (1,0.1) {\scriptsize $1$};
%\node[above] at (2,0.1) {\scriptsize $0$};

%\node[above] at (4,0.1) {\scriptsize $1$};
%\node[above] at (5,0.1) {\scriptsize $0$};
%\node[above] at (6,0.1) {\scriptsize $1$};
%\end{tikzpicture}

Множину всіх можливих префіксів позначимо $ P := \{0,1\}^M$. Звернемо увагу на те, що префіксів скінченна кількість. Тому для фіксованого рангу $n$ існує лише скінченна кількість попарно різних (з точністю до паралельного переносу) множин $\Delta_{\alpha_1\ldots\alpha_n}$, адже при однаковому префіксі множини є рівними з точністю до паралельного переносу. 
Для різних рангів $n$ та $m$ --- якщо  $\nabla_{\alpha_1\ldots\alpha_n}^s$ та
$\nabla_{\beta_1\ldots\beta_m}^s$
мають однакові префікси, то 
$\Delta_{\alpha_1\ldots\alpha_n}$ та $\Delta_{\beta_1\ldots\beta_m}$
є подібними з коефіцієнтом $s^{m-n}$.

\subsection*{Перехід від рангу $n$ до рангу $n+1$}

Опишемо перехід від циліндричних піввідрізків рангу $n$ (при поданні $C = E_n \oplus R_n$) до циліндричних піввідрізків рангу $n+1$ (при поданні $C = E_{n+1} \oplus R_{n+1}$).

%$l = \sum_{i=1}^{n} \frac{\alpha_i}{s^i} =\Delta_{\alpha_1 \ldots \alpha_n (0)}$
Нехай $(\alpha_1 \alpha_2 \ldots \alpha_n)$ --- фіксований набір елементів з алфавіту $A_s$, $l = \sum_{i=1}^{n} \frac{\alpha_i}{s^i}$. Тоді для довільного $d \in A_s$ маємо

\begin{equation}\label{transition_B}
B(\nabla_{\alpha_1 \alpha_2 \ldots \alpha_n d}^s) =
\Bigl(B(\nabla_{\alpha_1 \alpha_2 \ldots \alpha_n}^s) 
\oplus \tfrac{V}{s^{n+1}} \Bigr) 
\cap \Bigl[l +\tfrac{d-M+1}{s^{n+1}},l + \tfrac{d}{s^{n+1}}\Bigr]
\end{equation}

\begin{equation}\label{transition_SetOfSubsums}
\Delta_{\alpha_1 \ldots \alpha_n d} = 
\Bigl( B(\nabla_{\alpha_1 \ldots \alpha_n d}^s) \oplus R_{n+1} \Bigr) 
\cap \nabla_{\alpha_1 \ldots \alpha_n d}^s. 
\end{equation}

Префіксні точки $\{t_0,\ldots,t_{M-1}\}$ підциліндра $\nabla_{\alpha_1 \ldots \alpha_n d}$ лежать у відрізку $[l +\tfrac{d-M+1}{s^{n+1}},l + \tfrac{d}{s^{n+1}}]$, обидва кінці --- це самі префіксні точки ($t_0$ і $t_{M-1}$).
Префікси підциліндрів (прямих нащадків) однозначно визначаються через префікс їх прямого предка. Цей результат природний, адже префікс циліндра повністю визначає структуру $C$ на ньому, отже, і на всіх його підциліндрах.

Зазначимо, що наведені міркування про структуру множини $C$ зберігаються, якщо початковим (рангу $0$) циліндричним піввідрізком  розглядати не лише $[0,1)$, а й будь-який піввідрізок $[k,k+1), k \in \{0,1,\ldots,M-1\}$, описану схему можна використати для дослідження топологічної структури множини $C$ на кожному з них.

\section{Графи переходів}

\begin{definition}
Нехай $W$ --- непорожня множина. Відношення $R \subset W \times W$ називається \emph{відношенням часткової еквівалентності}~\cite{Mitchell_1996_PER}, якщо воно є симетричним $(aRb \;\Rightarrow\; bRa)$ і транзитивним $(aRb \ \text{і}\ bRc \;\Rightarrow\; aRc)$.

%тобто для всіх $a,b,c \in W$ виконуються умови:

\end{definition}
%При цьому рефлексивність, взагалі кажучи, не вимагається, тобто $aRa$ може не виконуватися для деяких $a$. 

Для елемента $a$ з $aRa$ \emph{класом часткової еквівалентності} назвемо множину ${b \in W : aRb}$.
Зокрема зазначимо, що на множині $W$ можуть існувати елементи без властивості $aRa$ (рефлексивність), тобто такі, що не входять до жодного класу часткової еквівалентності, індукованого $R$. 

%\subsection*{Відношення часткової еквівалентності для орієнтованого графа.}

Нехай $G$ --- орієнтований граф. Для вершин $a,b$ графа $G$ покладемо
$ a \sim b$, 
якщо існують орієнтовані шляхи з вершини $a$ до вершини $b$ та з вершини $b$ до вершини $a$. Інакше кажучи, $a \sim b$, якщо існує цикл (замкнений орієнтований шлях), що проходить і через $a$, і через $b$. 

Відношення $\sim$ є відношенням часткової еквівалентності на множині вершин графа, а його класи відповідають нетривіальним сильно зв'язним компонентам графа. Вершини, що не лежать на жодному циклі, не належать до жодного класу частокової еквівалентності. Назвемо їх \emph{транзитними} вершинами.

\begin{definition}\label{graph_of_transitions}
\emph{Графом переходів} $G$ називатимемо орієнтований граф з множиною вершин $P = \{0,1\}^M$ (вершини --- всі можливі префікси), у якому з кожної вершини $q \in P$ для кожного $d \in \{0,1,\ldots,s-1\}$ виходить ребро $q \rightarrowtail q_d$, де $q_d$ --- префікс $d$-го підциліндра довільного циліндричного піввідрізка з префіксом $q$.
\end{definition}

Переміщення по ребру орієнтованого графа --- це вибір відповідного підциліндра, тож рух по ребрах графа визначає послідовність вкладених підциліндрів. Символом $\rightarrowtail$ тут і надалі позначається окреме орієнтоване ребро. Означення коректне: префікс $q_d$ не залежить від вибору піввідрізка з префіксом $q$ (див. формулу~\eqref{transition_B}).

З кожної вершини графа $G$ виходить рівно $s$ ребер; допускаються петлі (префікс $d$-го підциліндра може збігатися з префіксом самого піввідрізка) та кратні ребра (різні значення $d$ можуть давати той самий префікс-нащадок).

Будемо писати $a \rightsquigarrow b$, якщо в графі $G$ існує орієнтований шлях з вершини $a$ до вершини $b$. 
Як описано вище, відношення $\sim$ ($a \sim b$, якщо $a \rightsquigarrow b$ та $b \rightsquigarrow a$) розбиває вершини графа на класи часткової еквівалентності та транзитні вершини. Стягнемо кожен клас часткової еквівалентності (тобто кожну нетривіальну сильно зв'язну компоненту графа) до однієї вершини; ця процедура називається конденсацією графа \cite{cormen2022introduction}. У результаті отримаємо орієнтований ациклічний граф $\bar{G}$, його вершинами є класи часткової еквівалентності та транзитні вершини. 

\begin{definition} 
Клас часткової еквівалентності називається \emph{термінальним}, якщо всі ребра, що виходять з його вершин, ведуть до вершин цього ж класу (інакше кажучи, з відповідної вершини конденсованого графа $\bar{G}$ не виходить жодне ребро).
\end{definition}

\begin{example}
Клас, що складається з вершини \((0,0,\ldots,0)\), є термінальним. 
Дійсно, префікс $(0,\ldots,0)$ означає, що $B(\nabla_{\alpha_1\ldots\alpha_n}^s)= \varnothing$. 
За~\eqref{transition_B} $B(\nabla_{\alpha_1\ldots\alpha_n d}^s) = \varnothing$ для всіх $d \in \{0,\ldots,s-1\}$, тобто префікси  підциліндрів також $(0,\ldots,0)$. Отже, всі $s$ ребер з цієї вершини є петлями: вершина $(0,\ldots,0)$ утворює клас, і жодне ребро не виходить за його межі.
\end{example}

\begin{definition}
\emph{Рівнем} класу часткової еквівалентності або транзитної вершини називається максимальна довжина орієнтованого шляху в конденсованому $\bar{G}$, що починається в даному класі/вершині та закінчується в деякому термінальному класі.
\end{definition}

\textbf{Зауваження.} Термінальні класи --- це класи рівня $0$.

\textbf{Зауваження.} Якщо з класу/транзитної вершини $X$ можна перейти до іншого класу/транзитної вершини $Y$, то рівень $X$ строго більший за рівень $Y$.

\begin{definition}
\emph{Стартовими префіксами} назвемо префікси циліндричних піввідрізків нульового рангу $[k-1, k)$, $k = 1, \ldots, M$, тобто набори
\[
(0,0,\ldots,0,1),\ (0,0,\ldots,1,0),
\ \ldots,\ (1,0,\ldots,0,0).
\]
В графі $G$ стартові префікси називатимемо \emph{стартовими вершинами}.
\end{definition}

За побудовою графа $G$ для кожного $n$ префікси циліндричних піввідрізків рангу $n$, що містяться в $[k-1, k)$, --- це в точності ті вершини, до яких у $G$ існує шлях довжини $n$ зі стартової вершини піввідрізка $[k-1, k)$. Тому для дослідження топологічної структури множини $C$ достатньо стартових вершин та вершин, досяжних з них.

\begin{definition}
Через $G'$ позначимо підграф графа $G$ з множиною вершин $P' \subseteq P$, що складається зі стартових вершин та вершин, досяжних з них орієнтованими шляхами. Ребрами $G'$ є всі ребра $G$, що виходять з вершин $P'$; через $\bar{G'}$ позначимо конденсацію графа $G'$.

Разом з кожною вершиною граф $G'$ містить усіх її нащадків, тому термінальність класів та рівні вершин у $\bar{G'}$ ті самі, що й у $\bar{G}$. Надалі працюємо з графами $G'$ та $\bar{G'}$.
\end{definition}

\section{Критерій наявності відрізка в $C$}

\begin{definition}
\emph{Прогалиною} множини $C$ називатимемо відкритий інтервал
$
(a,b)\subset [0,M], % можливо заміню на [0, r/(s-1)]
\qquad
(a,b)\cap C=\varnothing.
$
\end{definition}
\begin{remark}
В деяких статтях прогалиною називають максимальний відкритий інтервал у доповненні до $C$ ~\cite{Nowakowski_OnANewCondition}. Ми ж використовуємо термін у ширшому значенні "будь-який відкритий інтервал, що не перетинається з $C$".
\end{remark}

\begin{lemma}\label{lemma_gap}
Нехай $\nabla_{\alpha_1\ldots\alpha_n}^s \subset [0,\frac{r}{s-1}]$ має префікс $q$. Піввідрізок $\nabla_{\alpha_1\ldots\alpha_n}^s$ містить прогалину тоді й лише тоді, коли в графі префіксів $G'$ існує орієнтований шлях $ q \;\rightsquigarrow\; (0,0,\ldots,0).$
\end{lemma}

\begin{proof}

\emph{Достатність $(\Leftarrow)$.}
Нехай шлях $q \rightsquigarrow (0,\ldots,0)$ має довжину $k$. За побудовою графа цей шлях відповідає вибору цифр $d_1,\ldots,d_k$, тож циліндричний підінтервал $\nabla_{\alpha_1\ldots\alpha_n d_1\ldots d_k}^s$ має префікс $(0,\ldots,0)$. Тоді $B(\nabla_{\alpha_1\ldots\alpha_n d_1\ldots d_k}^s) = \varnothing$ і, відповідно, $\Delta_{\alpha_1\ldots\alpha_n d_1\ldots d_k} = \varnothing$.

\emph{Необхідність $(\Rightarrow)$.} Нехай в $\nabla_{\alpha_1\ldots\alpha_n}^s$ міститься прогалина $(a,b)$. Виберемо $N > n$ так, щоб виконувалась нерівність $\tfrac{1}{s^N} < \tfrac{b-a}{3}$. Для такого $N$ існує підциліндр
$\nabla_{\alpha_1\ldots\alpha_n\beta_{n+1}\ldots\beta_N}^s \subset (a,b)$. 
Виберемо число $k$ достатньо велике, щоб $s^k> M+1$. Далі ітеративно $k$ разів виберемо найправіший підциліндр у $\nabla_{\alpha_1\ldots\alpha_n\beta_{n+1}\ldots\beta_N}^s$, тобто $\nabla_{\alpha_1\ldots\alpha_n\beta_{n+1}\ldots\beta_N [s-1]\ldots[s-1]}^s$. Він і всі $M$ його префіксних точок лежать у $\nabla_{\alpha_1\ldots\alpha_n\beta_{n+1}\ldots\beta_N}^s$, отже, лежать у $(a,b)$. Тому його префікс $(0,0,\ldots,0)$. Отже, ми вказали орієнтований шлях від $q$ до $(0,0,\ldots,0)$ в графі $G'$.
\end{proof}

\begin{lemma}\label{lemma_black}
Нехай циліндричний інтервал
$\nabla_{\alpha_1\ldots\alpha_n}^s$
має префікс $q$.
Піввідрізок $ \nabla_{\alpha_1\ldots\alpha_n}^s \subset C$
тоді й лише тоді, коли
в графі префіксів $G'$ не існує орієнтованого шляху
$
q \rightsquigarrow (0,0,\ldots,0).
$
\end{lemma}

\begin{proof}

\emph{Необхідність $(\Rightarrow)$.} З існування шляху $q \rightsquigarrow (0,\ldots,0)$ випливає існування прогалини у $\nabla_{\alpha_1\ldots\alpha_n}^s$ (лема~\ref{lemma_gap}).

\emph{Достатність $(\Leftarrow)$.}
Припустимо, що орієнтованого шляху $q \rightsquigarrow (0,0,\ldots,0)$ не існує.
Візьмемо довільну точку $x_0 \in \nabla_{\alpha_1\ldots\alpha_n}^s$ та розглянемо послідовність циліндрів 
$ \nabla_{\alpha_1\ldots\alpha_n }^s \supset
 \nabla_{\alpha_1\ldots\alpha_n \alpha_{n+1}}^s   \supset 
 \nabla_{\alpha_1\ldots\alpha_n 
\alpha_{n+1} \alpha_{n+2}}^s \supset \ldots$,
де $\alpha_{n+1},\alpha_{n+2},\ldots$ --- наступні цифри деякого $s$-кового запису $x_0$. 

%Якщо $x_0 = \tfrac{r}{s-1}$, то $x_0 \in C$ (досягається при всіх $\alpha_i = r$), тому далі вважаємо, що $x_0 < \tfrac{r}{s-1}$.
%Якщо $x_0$ є $s$-ково-раціональним та $x_0 > 0$, тобто має вигляд $\frac{a}{s^j}, a > 0$, то воно має два різних зображення, а саме з періодами $(0)$ та $(s-1)$. Для однозначності обираємо його зображення з періодом $(0)$.

Вкладені циліндри --- це послідовність вкладених піввідрізків $[a_k,b_k)$. Якби у деякого $[a_k,b_k)$ префікс був $(0,0,\ldots,0)$, то існував би шлях $q \rightsquigarrow (0,\ldots,0)$, що дає суперечність з припущенням. Відповідно, у них ненульові префікси. Отже, кожен з відрізків $[a_k - \tfrac{M-1}{s^{n+k}}, a_k]$ містить як мінімум одну префіксну точку 
$x_k \in C$. Має місце нерівність
$|x_0 - x_k| 
\leq \tfrac{M-1}{s^{n+k}} + \tfrac{1}{s^{n+k}} =
\tfrac{M}{s^{n+k}}$,
тобто $x_k \to x_0$. Множина $C$ --- замкнена (схема доведення в ~\cite{Kakeya_1914_partial_sums}), тому і $x_0 \in C$. Точка $x_0$ обиралась довільно з $\nabla_{\alpha_1\ldots\alpha_n}^s$, отже, $\nabla_{\alpha_1\ldots\alpha_n}^s \subset C$.

\end{proof}

\begin{theorem}\label{theorem_segment_existence}
Множина $C$ містить відрізок тоді й лише тоді, коли конденсований граф $\bar{G'}$ містить термінальний клас
часткової еквівалентності, відмінний від класу, що складається з вершини
$(0,0,\ldots,0)$.
\end{theorem}
\begin{proof}
\emph{Необхідність $(\Rightarrow)$.} Нехай $C$ містить відрізок $[c,d]$. Виберемо досить велике $n$, щоб 
 $\tfrac{1}{s^n} < \tfrac{d-c}{3}$.
 Тоді існує $\nabla_{\alpha_1\ldots\alpha_n}^s \subset [c, d]$, позначимо його префікс $q$. Оскільки $\nabla_{\alpha_1\ldots\alpha_n}^s \subset C$, то він не містить прогалин, тож за лемою~\ref{lemma_gap} в графі $G'$ не існує шляху $q \rightsquigarrow (0,\ldots,0)$.
%Оскільки $q$ --- префікс циліндра, то $q \in P'$. 
Розглянемо клас/транзитну вершину у графі $\bar{G'}$, що містить $q$. Оскільки в конденсованому графі $\bar{G'}$ рівень при переходах завжди зменшується, то з $q$ існує шлях до термінального класу, котрий, як показано, не $\{(0,\ldots,0)\}$.

\emph{Достатність $(\Leftarrow)$.} Нехай $T$ --- термінальний клас у $\bar{G'}$, відмінний від $\{(0,\ldots,0)\}$, $q \in T$ --- довільна вершина. Оскільки $T$ термінальний, то не існує шляху $q \rightsquigarrow (0,0,\ldots,0)$. Оскільки $q \in P'$, то існує деякий циліндричний піввідрізок з префіксом $q$ і, згідно з попередньою лемою~\ref{lemma_black}, він повністю лежить у $C$. Всередині піввідрізка можна вибрати відрізок.
\end{proof}

%--------------------------------------------
\section{Критерій скінченності кількості відрізків у $C$}

\begin{definition}
Означимо \emph{кольори класів часткової еквівалентності та транзитних вершин} наступним чином:

\emph{Рівень $0$.} 
Термінальному класу $\{(0,0,\ldots,0)\}$ присвоюється \emph{білий} колір, усім іншим термінальним класам --- \emph{чорний}.

\emph{Рівні $\geq 1$.}
Класу (транзитній вершині) присвоюється колір залежно від кольорів безпосередніх нащадків, а саме:
\begin{itemize}
  \item \emph{білий}, якщо всі безпосередні нащадки мають білий колір;
  \item \emph{чорний}, якщо всі безпосередні нащадки мають чорний колір;
  \item \emph{чорно-білий} в інших випадках, тобто якщо серед безпосередніх нащадків є і чорний, і білий, або хоча б один чорно-білий нащадок.
\end{itemize}
\end{definition}
\begin{definition}
    \emph{Кольором вершини} з деякого класу еквівалентності покладемо колір класу еквівалентності, якому вона належить.
\end{definition}

\begin{remark}
    Інтуїтивно: білий колір позначає присутність прогалин у нащадках, а чорний --- присутність відрізків.
\end{remark}

\begin{remark}\label{remark_black_is_full}
Циліндричний інтервал $\nabla_X^s$ повністю лежить у $C$ тоді й лише тоді, коли префікс $X$ --- чорний (лема~\ref{lemma_black}).
\end{remark}

%Для префікса (вершини) $X$ через $\nabla_X^s$ позначатимемо довільний циліндричний піввідрізок із префіксом $X$. Оскільки множина $\Delta_{\alpha_1\ldots\alpha_n}$ з точністю до подібності визначається лише префіксом, то всі властивості, які нас цікавитимуть, не залежать від вибору такого циліндра.

\begin{remark}\label{remark_white_no_segment}
Циліндричний інтервал $\nabla_W^s$ з білим префіксом $W$ не містить відрізка
з $C$.

%Справді, якби існував відрізок $[a,b] \subset \nabla_W^s \cap C$, то всередині $[a,b]$ лежав би деякий циліндричний інтервал $\nabla_{q^*}^s$ достатньо великого рангу. Тоді $\nabla_{q^*}^s \subset C$. За попереднім зауваженням \ref{remark_black_is_full} $q^*$ --- чорний, тож існує шлях з білого $w$ до чорного $q^*$, отримали суперечність.

Також зазначимо, що білий колір префікса свідчить лише про відсутність відрізків у циліндричному піввідрізку, в той же час він може містити багато точок з $C$.
\end{remark}

\begin{lemma}\label{lemma_color_inheritance}
Колір будь-якої вершини містить кольори всіх своїх нащадків у такому сенсі:
якщо в графі $G'$ існує шлях між вершинами $q_1 \rightsquigarrow q_2$, то
\begin{itemize}
  \item якщо $q_2$ --- чорна, то $q_1$ --- чорна або чорно-біла;
  \item якщо $q_2$ --- біла, то $q_1$ --- біла або чорно-біла;
  \item якщо $q_2$ --- чорно-біла, то $q_1$ --- чорно-біла.
\end{itemize}
\end{lemma}

\begin{proof}
Нехай $q'$ --- передостання вершина шляху $q_1 \rightsquigarrow q_2$.
У $G'$ існує ребро $q' \rightarrowtail q_2$, тож або $q'$ і $q_2$ лежать в одному класі еквівалентності, й тому мають однаковий колір, або ж у $\bar{G'}$ клас/транзитна вершина, що містить $q_2$,
є безпосереднім нащадком класу/транзитної вершини, що містить $q'$, і тоді, за означенням кольору, колір вершини $q'$ містить колір вершини $q_2$. За індукцією по довжині шляху отримуємо, що колір вершини $q_1$ містить колір вершини $q_2$. 
\end{proof}

\subsection{Глобальна структура чорно-білих класів}

\begin{definition}
Клас часткової еквівалентності називається
\emph{циклічним}, якщо з кожної його вершини виходить рівно одне ребро, що веде до вершини цього ж класу (решта $s-1$ ребер виходять за межі класу).
\end{definition}

\begin{lemma}\label{infinite_paths}
Якщо $\gamma_i:a \rightsquigarrow b, \, a,b \in G', \, i \in \mathbb{N}$ такі шляхи, що будь-яка пара відрізняється на деякому кроці і вершина $b$ --- чорно-біла, то циліндричний піввідрізок $\nabla_{a}$ (з префіксом $a$) містить нескіченну кількість відрізків з $C$.
\end{lemma}

\begin{proof}
Оскільки $b$ --- чорно-біла, то з неї існує шлях до термінального класу $\{(0,\ldots,0)\}$ (білий нащадок) і шлях до іншого термінального класу (чорний нащадок). За лемами~\ref{lemma_gap} та~\ref{lemma_black}, циліндричний інтервал з префіксом $b$ містить як відрізок (у $C$), так і прогалину.

Шляхи $\gamma_i$ відрізняються на деякому кроці --- жоден не є початком іншого, тому циліндри-нащадки піввідрізка $\nabla_{a}$, що отримані переходами $\gamma_i$ є неперетинними. Для $\gamma_1, \ldots, \gamma_k$ розташуємо циліндри-нащадки на числовій прямій, в циліндрах на непарних місцях виділимо відрізок з $C$, а на парних --- прогалину. Оскільки $k$ довільне, то $\nabla_{a}$ містить нескінченну кількість відрізків з $C$.
\end{proof}

\begin{lemma}\label{necessary_condition1_branching}
Якщо в графі $G'$ існує нециклічний чорно-білий клас,
то множина $C$ містить нескінченну кількість відрізків.
\end{lemma}

\begin{proof}
Нехай $K \in \bar{G'}$ --- нециклічний чорно-білий клас, $x \in K$ ---довільна вершина.  
%і циліндричний піввідрізок $\nabla_x$.
За означенням нециклічності існує вершина $y \in K$, що має два внутрішні вихідні ребра $y \rightarrowtail a_1$, $y \rightarrowtail a_2$, де $a_1, a_2 \in K$. Клас $K$ є сильно зв'язною компонентою, тому існують шляхи $x \rightsquigarrow y$, $a_1 \rightsquigarrow x$ та $a_2 \rightsquigarrow x$, причому оберемо такі, що мають вершину $x$ лише на краях. Скомпонуємо їх таким чином:
\[\gamma_{a_1} = x \rightsquigarrow y \rightarrowtail a_1 \rightsquigarrow x; \quad
\gamma_{a_2} = x \rightsquigarrow y \rightarrowtail a_2 \rightsquigarrow x.
\]

\begin{center}
\begin{tikzpicture}[
  dot/.style={circle, fill, inner sep=2pt},
  cls/.style={draw, line width=1.2pt},
  arr/.style={->, line width=1.1pt}
]

\draw[cls] (0,0) ellipse (2.6 and 1.4);
\node at (3.5,0) {$K$};

\node[dot] (x1) at (-1.5,0) {};
\node[left] at (x1) {$x$};
\node[dot] (x2) at (0,-0) {};
\node[dot] (x3) at (1.5,0) {};
\node[right] at (x3) {$y$};

\node[dot] (x41) at (1,1) {};
\node[right] at (x41) {$a_1$};
\node[dot] (x42) at (0.75,-1) {};
\node[right] at (x42) {$a_2$};
\node[dot] (x52) at (-0.75,-1) {};

\draw[arr, blue] (x1) to[bend left=15] (x2);
\draw[arr, red]  (x1) to[bend right=15] (x2);
\draw[arr, blue] (x2) to[bend left=15] (x3);
\draw[arr, red]  (x2) to[bend right=15] (x3);

\draw[arr,blue] (x3) -- (x41);
\draw[arr,blue] (x41) -- (x1);

\draw[arr,red] (x3) -- (x42);
\draw[arr,red] (x42) -- (x52);
\draw[arr,red] (x52) -- (x1);

\node[blue] at (-0.1,0.8) {$\gamma_{1}$};
\node[red] at (-0.1,-0.8) {$\gamma_{2}$};

\end{tikzpicture}
\end{center}
Задамо $\gamma_k = (\gamma_{a_1})^k \gamma_{a_2} : x \rightsquigarrow x$. Будь-яка пара відрізняється на деякому кроці, вершина $x$ чорно-біла, тому за лемою~\ref{infinite_paths} $\nabla_{\alpha_1\ldots\alpha_n}^s$ містить нескінченну кількість відрізків з $C$.
\end{proof}

\begin{lemma}\label{necessary_condition2_connecting}
Якщо для двох різних чорно-білих класів з $G'$ існує орієнтований шлях від одного з них до іншого, то множина $C$ містить нескінченну кількість відрізків.
\end{lemma}
\begin{proof}
Нехай для різних чорно-білих класів $X, Y$ існує шлях $\gamma_{xy} \colon x \rightsquigarrow y$, де $x \in X$, $y \in Y$. Позначимо $\gamma'_{xy}$ деякий його підшлях, що містить $x$ та $y$ лише на краях. Зафіксуємо всередині класу $X$ деякий шлях $\gamma_{xx} \colon x \rightsquigarrow x$, що має $x$ лише на краях.

\begin{center}
\begin{tikzpicture}[
  dot/.style={circle, fill, inner sep=2pt},
  cls/.style={draw, line width=1.2pt},
  arr/.style={->, line width=1.1pt}
]

\draw[cls] (0,0) ellipse (1.6 and 0.9);
\node at (2.5,0) {$X$};

\node[dot] (x1) at (-0.7,0.2) {};
\node[left] at (x1) {$x$};
\node[dot] (x2) at (0.7,-0.5) {};
\node[dot] (x3) at (0.7,0.4) {};

\draw[arr,blue] (x3) -- (x1);
\draw[arr,blue] (x1) -- (x2);
\draw[arr,blue] (x2) -- (x3);
\node[blue] at (-0.1,0.6) {$\gamma_{xx}$};

\node[dot] (transit1) at (-0.7,-1.2) {};
\draw[arr,red] (x1) -- (transit1);

\draw[cls] (0,-2.4) ellipse (1.6 and 0.9);
\node at (2.5,-2.4) {$Y$};

\node[dot] (y) at (-0.7,-2.4) {};
\node[left] at (y) {$y$};
\draw[arr,red] (transit1) -- (y);
\node[red] at (-1.1,-1.2) {$\gamma'_{xy}$};

\end{tikzpicture}
\end{center}

Задамо $\gamma_k = (\gamma_{xx})^k \gamma'_{xy} : x \rightsquigarrow y$ та застосуємо лему~\ref{infinite_paths}. Маємо, що $C$ містить нескінченну кількість відрізків.
\end{proof}

\subsection{Локальна структура чорно-білих класів}

\begin{definition}
Нехай ребро $v \rightarrowtail w$ графа $G'$ відповідає переходу від деякого піввідрізка $\nabla_{\alpha_1\ldots\alpha_n}^s$ до $\nabla_{\alpha_1\ldots\alpha_n d}^s$, 
($d \in \{0, 1, \ldots, s-1\}$). Число $d$ називатимемо \emph{індексом} ребра, позначатимемо його $\mathrm{ind}(v, w)$.
\end{definition}

Нехай $v$ --- вершина деякого циклічного класу часткової
еквівалентності $X$, $k$ --- індекс єдиного ребра переходу по класу до $v' \in X$.
Для $v$ розщепимо множину інших безпосередніх нащадків на дві множини
\[
S^{-}(v) := \{\, w \mid v \rightarrowtail  w,\ \mathrm{ind}(v,w) < k\},\]
\[
S^{+}(v) := \{\, w \mid v \rightarrowtail  w,\ \mathrm{ind}(v,w) > k\}.\]

% добавити сюди малюнок.

\begin{lemma}\label{Each_S_is_1colored}
Для скінченності кількості відрізків у $C$ необхідно, щоб для кожної вершини $v \in K$, де K --- чорно-білий клас з $G'$, кожна з множин $S^{-}(v)$ та $S^{+}(v)$ була одноколірною, тобто складалась лише з чорних або лише з білих вершин.
\end{lemma}
\begin{proof}
За лемою \ref{necessary_condition1_branching} вважаємо всі чорно-білі класи циклічними. 

Припустимо, що в деякому чорно-білому класі $K$ існує вершина $a$, для якої множина $S^{-}(a)$ або $S^{+}(a)$ не одноколірна.
Нехай це множина $S^{-}(a)$ (випадок $S^+(a)$ симетричний). Чорно-білих нащадків немає (лема \ref{necessary_condition2_connecting}), тож існують ребра $B \colon a \rightarrowtail b$, $W \colon a \rightarrowtail w$, де $b$ --- чорна, а $w$ --- біла, $b, w \in S^{-}(a)$. 

Зафіксуємо у класі $K$ деякий шлях $\gamma_{aa} \colon a \rightsquigarrow a$. Розглянемо шляхи 
\[\gamma_i^B \coloneqq (\gamma_{aa})^i B \colon a \rightsquigarrow b; \quad 
\gamma_i^W \coloneqq (\gamma_{aa})^i W \colon a \rightsquigarrow w, \qquad i \in \mathbb{N}.\]
Жоден шлях не є початком іншого: коли коротший шлях останнім кроком виходить з класу $K$, то довший все ще залишається у цьому класі; при однаковій довжині шляхи різні.

Для деякого $\nabla_a$ (з префіксом $a$) розглянемо підциліндри, що отримуються у результаті переходів $\gamma_i^B$ та $\gamma_i^W$. 
%За припущенням множина $S^{-}(a)$ не є одноколірною. 
Циліндри $\gamma_m^B(\nabla_a)$ та $\gamma_m^W(\nabla_a)$ лежать лівіше за 
$(\gamma_{aa})^{m+1}(\nabla_a)$.
Якщо $\mathrm{ind}(a, b) < \mathrm{ind}(a, w)$, то неперетинні циліндри зліва направо мають порядок
\[ \gamma_1^B(\nabla_a),
\gamma_1^W(\nabla_a),
\gamma_{2}^B(\nabla_a),
\gamma_{2}^W(\nabla_a), \ldots, 
\gamma_i^B(\nabla_a),
\gamma_i^W(\nabla_a),\]
тому $\nabla_a$ містить щонайменше $i$ (тобто нескінченну кількість) відрізків з $C$ (по відрізку в кожному з $i$ циліндрів з чорними префіксами, лема~\ref{lemma_black}), відділених прогалинами (циліндри з білими префіксами, лема~\ref{lemma_gap}). 
%Оскільки $i$ --- довільне, то $\nabla_a$ містить нескінченну кількість відрізків з $C$.

Дані міркування залишаються в силі і для $\mathrm{ind}(a, b) > \mathrm{ind}(a, w)$ (лише міняються місцями $\gamma_i^B(\nabla_a)$ та $\gamma_i^W(\nabla_a)$). 
\end{proof}

\begin{definition}
Казатимемо, що вершина $v$ має \emph{властивість}: 
\begin{itemize}
\item 
\emph{$D_1$}, якщо всі вершини з $S^{-}(v)$ чорні, а всі вершини з $S^{+}(v)$ білі;

\item 
\emph{$D_2$}, якщо всі вершини з $S^{-}(v)$ білі, а всі вершини з $S^{+}(v)$ чорні.
    
\end{itemize}
\end{definition}

\begin{remark}\label{empty_s1_or_s2}
Якщо одна з множин $S^{-}(v), S^{+}(v)$ одноколірна, а інша порожня, то виконується одна з властивостей $D_1$, $D_2$.
\end{remark}

\begin{lemma}\label{all_v_are_D1_or_D2_type}
Для скінченності кількості відрізків у $C$ необхідно, щоб кожна вершина з чорно-білого класу мала властивість $D_1$ або $D_2$.
\end{lemma}
\begin{proof}
Вважаємо всі чорно-білі класи циклічними (за лемою~\ref{necessary_condition1_branching}).

Припустимо, що в деякому чорно-білому класі $K$ існує вершина $a$, для якої не виконуються ні $D_1$, ні $D_2$. У силу леми \ref{Each_S_is_1colored} $S^{-}(a)$ та $S^{+}(a)$ одноколірні, тож порушення умов $D_1$ та $D_2$ означає, що кольори $S^{-}(a)$ та $S^{+}(a)$ збігаються. 
Не порушуючи загальності, вважатимемо колір чорним. Оскільки клас $K$ чорно-білий і, за лемою \ref{necessary_condition2_connecting}, не має чорно-білих нащадків, то з нього виходить ребро до білої вершини. Позначимо це ребро $W \colon a_2 \rightarrowtail w$, де $a_2 \in K$, а $w$ --- біла, $w \notin K$. Зауважимо, що $a_2 \neq a$, адже в $a$ лише чорні нащадки. Зафіксуємо ребро
$B \colon a \rightarrowtail b$, де $b$ --- чорна вершина поза $K$, причому виберемо його так, щоб $b \in S^+ (a)$, якщо $w \in S^+ (a_2)$ та $b \in S^- (a)$, якщо $w \in S^- (a_2)$. Такий вибір можливий, тому що $S^+(a) \neq \varnothing$ та $S^-(a) \neq \varnothing$, інакше за зауваженням \ref{empty_s1_or_s2} виконувалась би одна з властивостей $D_1, D_2$.

Зафіксуємо всередині класу $K$ цикл з вершини $a$, що проходить через вершину $a_2$, тобто $\gamma_{aa} := \gamma_{a a_2}\,\gamma_{a_2 a}$. Розглянемо шляхи
\[ \gamma_i^B \coloneqq (\gamma_{aa})^i\, B \colon a \rightsquigarrow b, \qquad
  \gamma_i^W \coloneqq (\gamma_{aa})^i\,\gamma_{a a_2}\, W \colon a \rightsquigarrow w, \qquad i \in \mathbb{N}. \]

Далі схема доведення така сама, як у попередній лемі. Жоден шлях не є початком іншого.
Для деякого циліндричного інтервалу $\nabla_a$ розглянемо підциліндри, що отримуються у результаті переходів $\gamma_i^B$ та $\gamma_i^W$. Якщо $w \in S^- (a_2)$ та, відповідно, $b \in S^- (a)$, то неперетинні підциліндри розташовані у порядку
$ \gamma_1^B(\nabla_a),
\gamma_1^W(\nabla_a),
\gamma_{2}^B(\nabla_a),
\gamma_{2}^W(\nabla_a), \ldots.$
Якщо ж $w \in S^+ (a_2)$ та, відповідно, $b \in S^+ (a)$, то маємо порядок 
$ \gamma_i^W(\nabla_a),
\gamma_i^B(\nabla_a),
\gamma_{i-1}^W(\nabla_a),
\gamma_{i-1}^B(\nabla_a), \ldots$.

Отже, $\nabla_a$ містить нескінченну кількість відрізків з $C$.
\end{proof}

\begin{lemma}\label{necessary_condition3_middle}
Для скінченності кількості відрізків у $C$ необхідно, щоб у кожному чорно-білому класі часткової еквівалентності всі вершини мали одну й ту ж властивість: або всі мають властивість $D_1$, або всі --- $D_2$.
\end{lemma}
\begin{proof}
Вважаємо необхідні умови скінченності кількості відрізків (леми~\ref{necessary_condition1_branching},~\ref{necessary_condition2_connecting}~\ref{all_v_are_D1_or_D2_type}) виконаними: всі чорно-білі класи циклічні; чорно-білий клас не має чорно-білого нащадка; всі вершини чорно-білих класів мають одну з властивостей $D_1$ або $D_2$.

Припустимо, що в деякому чорно-білому класі $K$ існують як вершини з властивістю $D_1$, так і вершини з властивістю $D_2$. Покажемо, що тоді можна вибрати такі дві вершини $a_1$ з властивістю $D_1$ та $a_2$ з властивістю $D_2$, що виконується хоча б одна з умов:

\begin{itemize}
\item[(1)]
$S^-(a_1) \neq \varnothing$ та $S^-(a_2) \neq \varnothing$;

\item[(2)]
$S^+(a_1) \neq \varnothing$ та $S^+(a_2) \neq \varnothing$.
\end{itemize}

Розщепимо на випадки залежно від того, чи існує в $K$ вершина з властивістю $D_1$ з $S^+ \neq \varnothing$ і чи існує вершина з властивістю $D_1$ з $S^- \neq \varnothing$:

\begin{description}
\item[(А) (так, так).] Існує $a_1^+: S^+(a_1^+) \neq \varnothing$ з властивістю $D_1$ та існує $a_1^-: S^-(a_1^-) \neq \varnothing$ з властивістю $D_1$. Нехай $a_2$ --- деяка вершина з властивістю $D_2$. Якщо $S^-(a_2) \neq \varnothing$, то шукана пара $(a_1^-, a_2)$, для неї виконується умова (1). Якщо ж $S^-(a_2) = \varnothing$, то, оскільки $|S^-(a_2) \cup S^+(a_2)| = s - 1 > 0$, маємо $S^+(a_2) \neq \varnothing$. Шуканою парою буде $(a_1^+, a_2)$, для неї виконується умова (2).

\item[(Б) (так, ні).] Існує $a_1^+: S^+(a_1^+) \neq \varnothing$ з властивістю $D_1$ та не існує вершини $a_1^-: S^-(a_1^-) \neq \varnothing$ з властивістю $D_1$. Відсутність вершин вигляду $a_1^-$ означає, що не існує вершини з властивістю $D_1$ з чорним вихідним ребром. Клас $K$ чорно-білий, тож існує вершина $a_2$ з властивістю $D_2$, що має чорне вихідне ребро, тобто $S^+(a_2)\neq \varnothing$. Шуканою парою є $(a_1^+, a_2)$, виконується умова (2). 

\item[(В) (ні, так).] Аналогічно (Б).

\item[(Г) (ні, ні).] Тоді для будь-якої вершини $a_1$ з властивістю $D_1$ маємо $S^-(a_1) = S^+(a_1) = \varnothing$, але всі вершини мають $s$ вихідних ребер, де $s \geq 2$, тому $S^-(a_1) \cup S^+(a_1) \neq \varnothing$. Маємо суперечність, отже, цей випадок неможливий.
\end{description}

Отже, доведено, що завжди існує пара вершин $(a_1,a_2)$ така, що виконується (1) або (2). Нехай це буде (1), тобто $S^-(a_1) \neq \varnothing$ та $S^-(a_2) \neq \varnothing$.

Для вершини $a_1$ з властивістю $D_1$ позначимо ребро з індексом $0$ як $B \colon a_1 \rightarrowtail b$. Оскільки $S^-(a_1) \neq \varnothing$, то $b \in S^-(a_1)$. З властивості $D_1$ всі вершини з $S^-(a_1)$ чорні, тож $b$ --- чорна вершина поза $K$. Так само для вершини $a_2$ з властивістю $D_2$: позначимо ребро з індексом $0$ як $W \colon a_2 \rightarrowtail w$. Тоді $w \in S^-(a_2) \neq \varnothing$, тож $w$ --- біла вершина поза $K$. 

Клас $K$ --- сильно зв'язний, зафіксуємо всередині класу $K$ шляхи $\gamma_{a_1 a_2}$, $\gamma_{a_2 a_1}$та позначимо цикл 
$\gamma_{a_1 a_1} \coloneqq \gamma_{a_1 a_2} \gamma_{a_2 a_1}$.
Розглянемо шляхи:
\[ \gamma_i^B \coloneqq (\gamma_{a_1 a_1})^i\, B \colon a_1 \rightsquigarrow b, \qquad
  \gamma_i^W \coloneqq (\gamma_{a_1 a_1})^i\,\gamma_{a_1 a_2}\, W \colon a_1 \rightsquigarrow w, \qquad i \in \mathbb{N}. \]

Для деякого $\nabla_{a_1}$ розглянемо підциліндри, що отримуються у результаті переходів $\gamma_i^B$ та $\gamma_i^W$. Як і у попередніх лемах, жоден шлях не є початком іншого, тож дані підциліндри неперетинні.

Всередині підциліндра $(\gamma_{a_1 a_1})^i (\nabla_{a_1})$ ребро $B$ має індекс $0$ і виводить у найлівіший підциліндр, а перший крок $\gamma_{a_1 a_2}$ має індекс більший за $0$, тож виводить у підциліндр, що знаходиться правіше від $\gamma_i^B(\nabla_{a_1})$. 
Так само маємо, що всередині підциліндра $\bigl((\gamma_{a_1 a_1})^i\gamma_{a_1 a_2}\bigr)(\nabla_{a_1})$ ребро $W$ має індекс $0$ та виводить у найлівіший підциліндр, а перший крок $\gamma_{a_2 a_1}$ --- у підциліндр, що знаходиться правіше від $\gamma_i^W(\nabla_{a_1})$. Тому неперетинні підциліндри, при розгляді зліва направо, розташовані у порядку
$ \gamma_1^B(\nabla_{a_1}),
\gamma_1^W(\nabla_{a_1}),
\gamma_{2}^B(\nabla_{a_1}),
\gamma_{2}^W(\nabla_{a_1}), \ldots.$
Оскільки $i$ --- довільне, то $\nabla_{a_1}$ містить нескінченну кількість відрізків.

Випадок виконання умови (2) симетричний випадку виконанню умови (1), тому міркування у випадку (2)  аналогічні.
\end{proof}

\begin{lemma}\label{lemma_one_segment}
Нехай всі вершини циклічного чорно-білого класу $K$ мають одну й ту ж властивість $D_1$ або $D_2$. 
Тоді кожен циліндричний піввідрізок з префіксом $q \in K$ містить рівно один (максимальний у межах циліндра) відрізок/піввідрізок з $C$.
\end{lemma}

\begin{proof}
Доведення проведемо для випадку $D_1$. Випадок $D_2$ симетричний і доводиться аналогічно.

Для довільної вершини $v \in K$ через $k = k(v)$ позначимо індекс єдиного внутрішнього ребра з неї, нехай $v \rightarrowtail v'$, $v' \in K$. Властивість $D_1$ означає, що вершини в $S^-(v)$ чорні, в $S^+(v)$ --- білі. Тому циліндр $\nabla_{v}$ зліва направо поділяється на три послідовні частини:
\begin{itemize}
  \item підциліндри з $S^-(v)$ (чорні префікси) --- послідовні піввідрізки, що повністю лежать у $C$ (зауваження~\ref{remark_black_is_full}). Разом вони утворюють піввідрізок, що лежить у $C$;
  \item підциліндр $\nabla_{v'}$, він теж має префікс з $K$;
  \item підциліндри з $S^+(v)$ (білі префікси) --- послідовні піввідрізки, кожен з яких не містить відрізків з $C$ (зауваження~\ref{remark_white_no_segment}). Разом вони утворюють піввідрізок, що не містить відрізків з $C$.
\end{itemize}

% TikZ: суцільна чорна смуга ліворуч, стрілка ``вглиб $C$, пусто посередині, хвиляста смуга праворуч.

Покладемо $u_0 \coloneqq q$ та зафіксуємо циліндричний інтервал $\nabla_{u_0} = [L_0, R_0)$. Позначимо послідовність вершин, що отримуються ітеративними переходами по внутрішньому ребру, як $u_1, u_2, \ldots$, тобто $u_i \overset{k(u_i)}\rightarrowtail u_{i+1}, i \geq 0$.

Послідовності $u_0, u_1, \ldots ,u_i, \ldots$ з $C$ відповідає послідовність вкладених циліндрів $\nabla_{u_i} = [L_i, R_i)$, для яких, як показано вище, 
$[L_i, L_{i+1}) \subset C$, $[ R_{i+1}, R_i )$ не містить жодного відрізка з $C$.

Послідовність $(L_i)$ неспадна, $(R_i)$ незростаюча, $R_i - L_i \to 0 (i \to \infty)$. 
Отже, вони прямують до деякої спільної точки $p \in [L_0, R_0]$. Об'єднуючи лівіші й правіші частини, маємо $[L_0, p) \subset C$, тоді як $(p, R_0)$ не містить жодного відрізка з $C$.
Має місце $p \neq L_0$ та $p \neq R_0$, інакше клас $K$ був би не чорно-білим.
%Припустимо, що $L_0 = p$. Тоді $L_i = L_0$ для всіх $i$, всі індекси внутрішніх ребер дорівнюють $0$, тому всі зовнішні ребра для кожної з вершин $(u_i)$ лежать у $S^+(u_i)$. Клас $K$ циклічний, тому в $(u_i)$ входять всі вершини класу. Тож за $D_1$ всі зовнішні ребра з класу $K$ ведуть до білих вершин. Таким чином з класу $K$ не існує шляху до чорної вершини, а це суперечить чорно-білості класу $K$. Аналогічно, припущення $p = R_0$ виключило б шлях до білої вершини. Отже, $L_0 < p < R_0$.
Отже, $\nabla_q$ містить рівно один максимальний у межах циліндра відрізок $[L_0, p]$ з $C$ ($p \in C$, бо $C$ --- замкнена).
\end{proof}

\begin{theorem}[Критерій скінченності кількості відрізків у $C$]\label{theorem_finiteness_criterion}
Множина $C$ містить лише скінченну кількість відрізків тоді і лише тоді, коли одночасно виконуються такі умови:

1) усі чорно-білі класи циклічні,

2) між чорно-білими класами немає шляхів,

3) для чорно-білого класу всі вершини в ньому мають одну й ту ж властивість: або всі мають властивість $D_1$, або всі --- $D_2$.
\end{theorem}

\begin{proof}
\emph{Необхідність.} Доведено в трьох лемах раніше (лема~\ref{necessary_condition1_branching}, лема~\ref{necessary_condition2_connecting}, лема~\ref{necessary_condition3_middle}).

\emph{Достатність.} Нехай умови (1--3) виконуються.

Казатимемо, що шлях $a \rightsquigarrow b$ є шляхом до деякого класу $K \ni b $, якщо $b$ є його першим (єдиним) потраплянням до класу $K$.

Припустимо, що шлях $p \rightsquigarrow K_1$ від деякої стартової вершини $p$ до чорно-білого класу $K_1$ проходить через деякий інший клас $K_0$, тобто $p \rightsquigarrow K_1 = p \rightsquigarrow K_0 \rightsquigarrow K_1$. Оскільки $K_1$ --- чорно-білий, то, за лемою~\ref{lemma_color_inheritance}, клас $K_0$ також чорно-білий. Це суперечить умові (2). Отже, шляхи від стартових вершин до чорно-білих класів проходять лише через транзитні вершини. 
Транзитні вершини не належать до жодного класу, отже, не належать до жодного циклу. Граф $G'$ має скінченну кількість транзитних вершин (взагалі вершин скінченна кількість). Чорно-білих класів та стартових вершин також скінченна кількість, тому існує лише скінченна кількість шляхів від стартових вершин до чорно-білих класів.

Кожен шлях $p \rightsquigarrow K$, де $p$ --- стартова вершина, $K$ --- чорно-білий клас, задає певний циліндричний піввідрізок $\nabla_X$, який за лемою~\ref{lemma_one_segment} містить рівно один відрізок з $C$ у своїх межах. 
Аналогічно, кількість шляхів $p \rightsquigarrow B$, де $p$ --- стартова вершина, $B$ --- чорний клас, що проходять лише через транзитні вершини у графі $G'$, також скінченна. Кожен шлях $p \rightsquigarrow B$ задає певний циліндричний піввідрізок $\nabla_Y \subset C$ (зауваження~\ref{remark_black_is_full}).

Описано два сімейства циліндричних піввідрізків, кожен з яких містить рівно один (пів)відрізок з $C$ у своїх межах. Їх скінченна кількість, тому об'єднання всіх цих циліндрів, котре позначимо $F$, містить лише скінченну кількість відрізків з $C$.

Припустимо, що існує відрізок $I \subset C$, що $I \cap F = \varnothing$. Виберемо всередині $I$ деякий циліндр $\nabla_u$; $\nabla_u \subset I \subset C$, тому префікс $u$ --- чорний (зауваження~\ref{remark_black_is_full}). 
Розглянемо шлях $p \rightsquigarrow u$.

Нехай шлях $p \rightsquigarrow u$ складається не лише з транзитних точок. Першу вершину на ньому, що належить деякому класу часткової еквівалентності, позначимо $u_1$. Вона має нащадка $u$, тому за лемою~\ref{lemma_color_inheritance} вона чорна або чорно-біла, тож потрапляє в одне з сімейств, описаних вище. 
Циліндричний піввідрізок, що відповідає шляху $p \rightsquigarrow u_1$, містить в собі $\nabla_u$, суперечність з припущенням.

Нехай шлях $p \rightsquigarrow u$ складається лише з транзитних вершин. Оскільки префікс $u$ чорний, то існує шлях $u \rightsquigarrow t$, де $t$ --- термінальний чорний клас. Позначимо $u_1$ першу вершину на ньому, що належить деякому класу часткової еквівалентності. Оскільки $u$ чорний, то всі його нащадки чорні, отже, $u_1$ --- чорний. 
Шлях $p \rightsquigarrow u\rightsquigarrow u_1$ задає циліндричний інтервал, що є елементом другого сімейства, отримано суперечність.

Таким чином описані два сімейства вичерпують всі відрізки з $C$.

Отже, $C$ містить лише скінченну кількість відрізків.
\end{proof}

\section{Приклади застосування}

Для мультигеометричних рядів $(k_1, k_2, \ldots, k_m; \frac{1}{s})$ критерій існування відрізка та критерій скінченності їх кількості в множині підсум дають можливість скористатись класифікаційною теоремою Гатрі-Німана та дізнатись топологічну структуру множини підсум, тобто для мультигеометричних рядів вигляду~\eqref{multigeom_series} отримано критерій канторвальності.

Не кожна множина $C$ є множиною неповних сум ряду, 
тому теорема Гатрі-Німана стосовно неї безпосередньо не застосовна. %, тобто з умови існування нескінченної кількості відрізків канторвальність $C$ не випливає.

\emph{Приклад 1}. Ряд Гатрі-Німана: $A = \{0,2,3,5\}, s=4$. Граф $G$:

{\centering
\begin{tikzpicture}[
  vertex/.style={circle, draw, minimum size=10pt},
  arr/.style={->, line width=1pt}
]

% ---- вершини
\node[vertex] (v00) at (-4,0) {$(0,0)$};
\node[vertex] (v10) at (-1.5,0) {$(1,0)$};
\node[vertex] (v01) at (1.5,0) {$(0,1)$};
\node[vertex] (v11) at (4,0) {$(1,1)$};

% (0,0) 
\draw[arr,->] (v00) to[out=60, in=120, looseness=3]   (v00);
\draw[arr,->] (v00) to[out=55, in=125, looseness=4]   (v00);
\draw[arr,->] (v00) to[out=50, in=130, looseness=5]   (v00);
\draw[arr,->] (v00) to[out=45, in=135, looseness=6]   (v00);

% (0,1) 
\draw[arr,->] (v01) to[out=60, in=120, looseness=3]   (v01);
\draw[arr,->] (v01) to[out=45, in=135, looseness=4]   (v01);
\draw[arr, bend left=15] (v01) to (v10);
\draw[arr, bend left=15] (v01) to (v11);

% (1,0) 
\draw[arr,->] (v10) to[out=60, in=120, looseness=3]   (v10);
\draw[arr,->] (v10) to[out=45, in=135, looseness=4]   (v10);
\draw[arr, bend left=15] (v10) to (v01);
\draw[arr, bend left=15] (v10) to (v00);

% (1,1) 
\draw[arr,->] (v11) to[out=60, in=120, looseness=3]   (v11);
\draw[arr,->] (v11) to[out=55, in=125, looseness=4]   (v11);
\draw[arr,->] (v11) to[out=50, in=130, looseness=5]   (v11);
\draw[arr,->] (v11) to[out=45, in=135, looseness=6]   (v11);

\end{tikzpicture}
\par}

Білий клас $(0,0)$, один чорний клас $(1,1)$.
Стартові вершини: $(0,1)$ для відрізка $[0,1)$ та $(1,0)$ для $[1,2)$, вони формують чорно-білий клас $\{ (0,1), (1,0)\}$. Існує шлях із стартової вершини до чорного класу, отже, множина підсум містить відрізок. Чорно-білий клас не є циклічним, отже, відрізків нескінченна кількість. За класифікаційною теоремою Гатрі-Німана множина підсум даного ряду --- канторвал.

\emph{Приклад 2}. Наведемо новий мультигеометричний ряд, канторвальність множини підсум якого не перевіряється вже відомими достатніми умовами, котрі описано в роботах~\cite{Bartoszewicz_Filipczak_Szymonik_2014} та ~\cite{Nowakowski_OnANewCondition}. Якщо обмежитись рядами вигляду (\ref{multigeom_series}), то ці достатні умови при переході до графа означають попадання в термінальний чорний клас $(1,1,\ldots,1)$. У наведеного далі ряду вершина $(1,1,\ldots,1)$ не є досяжною із стартових префіксів.

Ряд $(1,4,16; \tfrac{1}{5})$: $V = \{0,1,4,5,16,17,20,21\}$, $M = 6$. Конденсований граф $\bar{G'}$:

{\centering
\begin{tikzpicture}[
  scale=1.0,
  vertex/.style={ellipse, draw, inner sep=2pt, font=\scriptsize},
  cls/.style={draw, dashed, line width=0.9pt, rounded corners=9pt, inner sep=5pt, align=center, font=\scriptsize},
  arr/.style={->, line width=1pt}
]

% ---- стартові транзитні вершини
\node[vertex] (T1) at (-5.1,4.3) {\texttt{000100}};
\node[vertex] (T2) at (-2.7,4.3) {\texttt{001000}};
\node[vertex] (T3) at (0.9,4.3)  {\texttt{000010}};
\node[vertex] (T4) at (3.7,4.3)  {\texttt{010000}};

\node at (-5.1,4.8) {\scriptsize $[2,3)$};
\node at (-2.7,4.8) {\scriptsize $[3,4)$};
\node at (0.9,4.8)  {\scriptsize $[1,2)$};
\node at (3.7,4.8)  {\scriptsize $[4,5)$};

% чорно-білий клас 
\node[cls] (C1) at (-0.7,3.0) {\itshape чорно-білий клас\\[1pt]
\texttt{000001}\ \texttt{001100}\ \texttt{000110}\\
\texttt{100110}\ \texttt{110000}\ \texttt{100000}};

% транзитні вершини
\node[vertex] (t5) at (-3.4,2.0) {\texttt{000011}};
\node[vertex] (t6) at (2.6,2.0)  {\texttt{011000}};

% чорний клас
\node[cls] (K1) at (5,1.5) {\itshape чорний клас\\[1pt] \texttt{011001}};

% термінальні класи
\node[cls] (W) at (-4.8,-0.8) {\itshape білий термінальний клас\\[1pt] \texttt{000000}};

\node[cls] (B) at (1.0,-0.7) {\itshape чорний термінальний клас\\[1pt]
\texttt{110011}\ \texttt{100111}\ \texttt{001110}\ \texttt{011100}\\
\texttt{111001}\ \texttt{011101}\ \texttt{111011}\ \texttt{110111}\\
\texttt{101110}\ \texttt{001111}\ \texttt{011110}\ \texttt{111100}\\
\texttt{101111}\ \texttt{111101}};

% ребра 
% 000100
\draw[arr] (T1) to (C1);
\draw[arr] (T1) to (W);
% 001000
\draw[arr] (T2) to (C1);
\draw[arr] (T2) to (t5);
\draw[arr] (T2) to[bend right=25] (W);
% 000010
\draw[arr] (T3) to (C1);
\draw[arr] (T3) to (t6);
\draw[arr] (T3) to[bend left=20] (B);
% 010000
\draw[arr] (T4) to (C1);
\draw[arr] (T4) to (t6);
\draw[arr] (T4) to (K1);
\draw[arr] (T4) to[bend left=20] (B);
% чорно-білий клас
\draw[arr] (C1) to (t5);
\draw[arr] (C1) to (t6);
\draw[arr] (C1) to (W);
\draw[arr] (C1) to[bend right=15] (K1);
\draw[arr] (C1) to (B);
% транзитні вершини
\draw[arr] (t5) to (B);
\draw[arr] (t6) to (K1);
\draw[arr] (t6) to (B);
% чорний клас
\draw[arr] (K1) to (B);

\end{tikzpicture}
\par}

З стартових вершин досяжний чорний термінальний клас. При цьому чорно-білий клас не є циклічним: наприклад, з вершини $001100$ виходять чотири внутрішні ребра. 
Отже, множина підсум ряду містить нескінченну кількість відрізків, тобто є канторвалом.

\bibliographystyle{plain}
\bibliography{mybiblio-ed2}